\documentclass[10pt]{amsart}

\usepackage{amsthm,amsmath,amssymb,tikz}
\usetikzlibrary{matrix}
\usepackage[margin=3cm]{geometry}
\usepackage{amsmath,amssymb,amsthm,amsfonts,bbding,stmaryrd,dsfont,wasysym,pifont,xspace,xcolor}
\usepackage{parskip,fancyhdr,url}

\title{Shadows of Teichm\"uller discs in the curve graph}

\author{Robert Tang}
\address{Department of Mathematics, University of Oklahoma, Norman, OK 73019, USA}
\email{rtang@math.ou.edu}
\urladdr{http://www.math.ou.edu/~rtang/}

\author{Richard C. H. Webb}
\address{Department of Mathematics, University College London, Gower Street, WC1E 6BT, London, UK}
\email{richard.webb@ucl.ac.uk}
\urladdr{http://www.ucl.ac.uk/~ucahrch/}

\numberwithin{equation}{section}

\newtheorem{Thm}[equation]{Theorem}
\newtheorem{Prop}[equation]{Proposition}
\newtheorem{Lem}[equation]{Lemma}
\newtheorem{Cor}[equation]{Corollary}

\newtheorem{Que}{Question}

\theoremstyle{definition}
\newtheorem{Dfn}[equation]{Definition}

\theoremstyle{remark}
\newtheorem{Rem}[equation]{Remark}

\newcommand{\R}{\mathbb{R}}
\newcommand{\Z}{\mathbb{Z}}

\newcommand{\Hy}{\mathbb{H}}
\newcommand{\CC}{\mathbb{C}}

\newcommand{\diam}{\mathrm{diam}}
\newcommand{\area}{\mathrm{area}}
\newcommand{\uvec}{\mathbf{u}}

\newcommand{\width}{\mathrm{width}}
\newcommand{\height}{\mathrm{height}}
\newcommand{\length}{\mathrm{length}}
\newcommand{\sys}{\mathsf{sys}}
\newcommand{\cyl}{\mathsf{cyl}}
\newcommand{\hcyl}{\widehat{\mathsf{cyl}}}
\newcommand{\disc}{\Delta}

\newcommand{\C}{\mathcal{C}}
\newcommand{\X}{\mathcal{X}}
\newcommand{\Y}{\mathcal{Y}}
\newcommand{\F}{\mathcal{F}}
\newcommand{\M}{\mathcal{M}}
\newcommand{\MF}{\mathcal{MF}}
\newcommand{\PMF}{\mathcal{PMF}}

\newcommand{\B}{\mathcal{B}}
\newcommand{\MA}{\mathcal{MA}}
\newcommand{\FMA}{\mathcal{FMA}}
\newcommand{\A}{\mathcal{A}}
\newcommand{\arcd}{\mathcal{A}_\disc}

\newcommand{\mcg}{\mathrm{Mod}}

\newcommand{\bM}{\mathbf{M}}

\newcommand{\sK}{\mathsf{K}}
\newcommand{\sC}{\mathsf{C}}

\newcommand{\sQ}{\mathsf{Q}}

\newcommand{\sP}{\mathsf{P}}
\newcommand{\sL}{\mathsf{L}}
\newcommand{\sR}{\mathsf{R}}

\newcommand{\sh}{\mathsf{h}}
\newcommand{\sB}{\mathsf{B}}
\newcommand{\sA}{\mathsf{A}}
\newcommand{\sN}{\mathsf{N}}

\newcommand{\G}{\mathcal{G}}

\newcommand{\Teich}{\mathrm{Teich}}

\newcommand{\QD}{\mathcal{QD}}
\newcommand{\SL}{SL}
\newcommand{\SO}{SO}
\newcommand{\Ext}{\mathsf{Ext}}
\newcommand{\Hyp}{\mathsf{Hyp}}
\renewcommand{\L}{\mathcal{L}}

\newcommand{\hyp}{{\mathsf{\delta}}}
\newcommand{\widthB}{{\sf W}} 
\newcommand{\widthV}{{\sf W'}}
\newcommand{\projS}{{\sf P_{\sf 2}}}
\newcommand{\qcS}{{\sf Q_{\sf 2}}}
\newcommand{\lipS}{{\sf K_{2}}}
\newcommand{\lipcS}{{\sf C_{2}}}
\newcommand{\diamL}{{\sf D}}
\newcommand{\hdSC}{{\sf h_{3}}}
\newcommand{\hdsys}{{\sf h_{2}}}
\newcommand{\Gen}{{\sf G}}
\newcommand{\lipV}{{\sf K_{1}}}
\newcommand{\qcV}{{\sf Q_{1}}}
\newcommand{\projV}{{\sf P_{1}}}

\newcommand{\TT}{\mathcal{TT}}
\newcommand{\VD}{V_{\disc}}
\newcommand{\tq}{\tau_q}
\newcommand{\td}{\tau_\disc}

\newcommand{\proofof}[1]{\hfill\newline\noindent\emph{Proof of {#1}.} }
\newcommand{\halmos}{\hfill$\square$}

\begin{document}

	\begin{abstract}
	We consider several natural sets of curves associated to a given Teichm\"uller disc, such as the systole set or cylinder set, and study their coarse geometry inside the curve graph. We prove that these sets are quasiconvex and agree up to uniformly bounded Hausdorff distance. 	
	We describe two operations on curves and show that they approximate nearest point projections to their respective targets.
	Our techniques can be used to prove a bounded geodesic image theorem for a natural map from the curve graph to the filling multi-arc graph associated to a Teichm\"uller disc.
		\end{abstract}
	
	\maketitle
	
	\section{Overview}
	
	In their groundbreaking paper \cite{MM1}, Masur and Minsky proved that the curve graph is hyperbolic by studying its interplay with the large-scale geometry of Teichm\"uller space. Their notion of \emph{balance time} for curves on Teichm\"uller geodesics, in particular, proved useful for showing that Teichm\"uller geodesics ``shadow'' reparameterised quasigeodesics in the curve graph under the systole map. 
		
	Our paper studies the coarse geometry of \emph{Teichm\"uller discs}, a natural generalisation of Teichm\"uller geodesics, via various notions of ``shadows'' in the curve graph. We also generalise the balance time to \emph{balance points} for curves on Teichm\"uller discs, and use it to extend results of \cite{MM1}.
	
	Throughout this paper, let $S$ be a closed, connected, orientable surface of genus $g \geq 2$. The \emph{curve graph} $\C(S)$ associated to $S$ has as its vertices the free homotopy classes of non-trivial simple closed curves on $S$, with edges spanning vertices if the corresponding curves can be realised disjointly. We shall equip the \emph{Teichm\"uller space} $\Teich(S)$, the parameter space of marked conformal structures on $S$, with the \emph{Teichm\"uller metric}. A \emph{Teichm\"uller disc} $\disc = \disc_q$ is a geodesically embedded copy of the hyperbolic plane $\Hy^2$ (with curvature $-4$) in $\Teich(S)$ which parameterises the conformal structures in the $\SL(2,\R)$--orbit of a holomorphic quadratic differential $q$ on $S$. We mainly focus on the \emph{half-translation structures} on $S$ naturally associated to such quadratic differentials. 
	
	We shall use the adjective \emph{uniform} to describe constants depending only on (the topological type of) $S$, and \emph{universal} if they can be chosen independently of $S$.
	
	\begin{Thm}\label{shadows}
	    For any Teichm\"uller disc $\disc \subset \Teich(S)$, the following sets agree up to universal Hausdorff distance in $\C(S)$.
	    \begin{itemize}
		\item $V(\disc)$ -- the set of straight vertex cycles for $\disc$,
		\item $\sys(\disc)$ -- the flat systoles appearing on $\disc$,
		\item $\sys^\Ext(\disc)$ -- the extremal length systoles appearing on $\disc$,
		\item $\sys^\Hyp(\disc)$ -- the hyperbolic length systoles appearing on $\disc$.
	    \end{itemize}
	    The following agree with the aforementioned sets up to uniform Hausdorff distance.
	    \begin{itemize}
		\item $\cyl(\disc)$ -- the cylinder curves on some (hence every) $q \in \disc$,
		\item $\hcyl(\disc)$ -- the curves with constant direction on some (hence every) $q \in \disc$.
	    \end{itemize}
	\end{Thm}
	
	Most of the above sets have been studied in some form in the literature. The only new definition, to our knowledge, are straight vertex cycles. Informally speaking, $V(\disc)$ is the set of curves whose geodesic representatives on every $q \in \disc$ run over each saddle connection at most once in each direction and cannot be further decomposed via curve surgery. It is worth emphasising that curves in $V(\disc)$ can have arbitrarily large lower bounds on their flat length (and therefore the extremal and hyperbolic lengths) over $\disc$ -- in particular, they are not systoles. 	We provide a brief description, and refer the reader to Section \ref{secstraight} for more details. See also Section \ref{sectrain} for background on train tracks.
	
	A geodesic representative of a curve $\alpha \in \C(S)$ on a half-translation surface $q \in \disc$ is typically a concatenation of saddle connections. One can ``smooth'' such a geodesic representative at each singular point to obtain a train track $\tq(\alpha)$ carrying $\alpha$. This construction commutes with $\SL(2,\R)$--deformations of $q$, and so we can canonically define a train track $\td(\alpha)$. The set of \emph{straight vertex cycles} for $\disc$ is
	\[V(\disc) = \cup_{\alpha \in \C(S)} V(\td(\alpha)),\]
	where $V(\td(\alpha))$ denotes the vertex cycles of $\td(\alpha)$. 
	
	It is worth emphasising that this construction needs only the combinatorial pattern of the saddle connections used by geodesics, not their lengths or directions. Thus we can use combinatorial methods to show the following, with effective control on the constants.
			       
       \begin{Thm}\label{vertexset}
        There is a universal constant $\qcV$ so that for any Teichm\"uller disc $\disc$, the set $V(\disc)$ is $\qcV$--quasiconvex in $\C(S)$. Moreover, the operation $\alpha \mapsto V(\td(\alpha))$ agrees with the nearest point projection relation from $\C(S)$ to $V(\disc)$ up to universally bounded error. 
       \end{Thm}
       
       \begin{Cor}
        For any Teichm\"uller disc $\disc$, the systole sets $\sys(\Delta)$, $\sys^\Ext(\Delta)$, and $\sys^\Hyp(\Delta)$ are universally quasiconvex. The sets $\cyl(\disc)$ and $\hcyl(\disc)$ are uniformly quasiconvex. \halmos
       \end{Cor}

	Given $q \in \disc$, let $S' = S'(\disc)$ be the underlying topological surface with the singularities of $q$ considered as marked points.
			Any non-cylinder curve $\alpha \in \C(S)$ determines a multi-arc $\arcd(\alpha)$ on $S'$ by taking the set of saddle connections used by the geodesic representative of $\alpha$ on any $q\in \disc$. In particular, if $\alpha$ is sufficiently far from $V(\disc)$ in $\C(S)$ then $\arcd(\alpha)$ fills $S'$. Let $\FMA(S')$ be the graph whose vertices are filling multi-arcs on $S'$, with edges connecting disjoint pairs of filling multi-arcs, and $\FMA(\disc)$ its subgraph spanned by the multi-arcs with component arcs realisable as saddle connections on $q$. There is a natural $1$--Lipschitz retract $\FMA(S') \rightarrow \FMA(\disc)$ which maps a multi-arc to the set of saddle connections appearing on its geodesic representative on $q$. It follows that $\FMA(\disc)$ isometrically embeds into $\FMA(S')$. We prove the following bounded geodesic image theorem.

	  \begin{Thm} \label{introbgit}
	   There are universal constants $\sA$ and $\sB$ so that if $G$ is a geodesic in $\C(S)$ disjoint from the $\sA$--neighbourhood of $V(\disc)$, then $\arcd(G)= \{ \arcd(\alpha) ~|~ \alpha\in G\}$ has diameter at most $\sB$ in $\mathcal{FMA}(\Delta)$.
	  \end{Thm}

	The endpoints of a Teichm\"uller geodesic $\G$ correspond to a pair of measured transverse foliations which give the horizontal and vertical directions for every half-translation structure along $\G$.
	If a curve $\alpha$ in neither completely horizontal nor completely vertical along $\G$, then its \emph{balance time} on $\G$ is the unique point on $\G$ where the horizontal and vertical lengths of $\alpha$ are equal. Masur and Minsky show that the operation of taking a systole at the balance time of a curve satisfies certain coarse retraction properties \cite{MM1}.
	 
	We introduce the notion of a \emph{balance point} of a curve $\alpha$ on a Teichm\"uller disc $\disc$ in Section \ref{secbalance}. This is a half-translation structure $X \in \disc$ so that, under all possible rotations, the horizontal and vertical lengths of $\alpha$ agree up to a bounded ratio. Such a point can be used to approximate the balance times of $\alpha$ along \emph{all} Teichm\"uller geodesics on $\disc$.
	
	\begin{Thm}\label{thmbalance}
	Let $\disc$ be a Teichm\"uller disc and $\alpha \in \C(S)$ a curve. If $\alpha$ does not have constant direction on $\disc$ then there is a point $X \in \disc$, called a \emph{balance point} of $\alpha$ on $\disc$, satisfying the following.
	\begin{itemize}
	 \item For any Teichm\"uller geodesic $\G$ on $\disc $, the nearest point projection of $X$ to $\G$ in $\disc$ is at most a distance of $\log 2$ from the balance time of $\alpha$ on $\G$.
	 \item Any flat systole on $X$ is universally close to any nearest point projection of $\alpha$ to $\sys(\disc)$ in $\C(S)$.
	 \item If the length of $\alpha\in\C(S)$ is minimised at $X'\in\disc$, among all flat surfaces in $\disc$, then $X$ is universally close to $X'$.
	\end{itemize}
	\end{Thm}
	
	We also introduce the \emph{auxiliary polygon} $P_q(\alpha)$, associated to a curve $\alpha$ on a half-translation surface $q$, as a useful tool in proving this theorem by reducing our arguments to elementary Euclidean geometry in the plane. For example, the auxiliary polygon can be used to explicitly find a balance point on $\disc = \disc_q$ given a geodesic representative of $\alpha$ on $q$. The area of $P_q(\alpha)$ also gives an estimate for the minimal flat length of $\alpha$ on $\disc$.
	Further applications of auxiliary polygons are currently being investigated in a joint project involving the first author, Max Forester, and Jing Tao.
	
	\subsection*{Organisation}
	
	In Section \ref{secbackground}, we provide background on coarse geometry, curve graphs, train tracks, half--translation surfaces and Teichm\"uller discs. We then describe the construction of the train track $\tq(\alpha)$ in Section \ref{secstraight}, where we also give a proof of Theorem \ref{vertexset} and Theorem \ref{introbgit}.
   
	Theorem \ref{shadows} shall be proved through the following collection of coarse inclusions. In the diagram below, an arrow 
	$A \xrightarrow{r} B$
	indicates that $A$ is contained in the $r$--neighbourhood of $B$ in $\C(S)$, with arrows in both directions indicating coarse inclusions in both directions.
	\newline
	
	\begin{center}
	 \begin{tikzpicture}
	  \matrix (m) [matrix of math nodes,row sep=4em,column sep=6em,minimum width=2em]
	  {
	    \sys^\Ext(\disc) &  & V(\disc) & \\
	    \sys^\Hyp(\disc) & \sys(\disc) & \cyl(\disc) & \hcyl(\disc)\\};
		  \path[-stealth]
	    (m-1-1) edge [<->] node [left] {$37$} (m-2-1)
		    edge [<->] node [above] {$37$} (m-2-2)
	    (m-2-1.east|-m-2-2) edge [<->] node [below] {$37$}
		    (m-2-2.west)
	    				    	    (m-2-2) edge [dashed] node [above] {$\hdSC(g)$} (m-2-3)
		    edge [<-, bend right, dashed] node [below] {$1$} (m-2-3)
		    edge [double] node [right] {~$14$} (m-1-3)
		    edge [<-, bend left, double] node [above] {$54$} (m-1-3)
	    (m-2-3) edge [dotted, thick] node [above] {$0$} (m-2-4)
		    edge [<-, bend right, dashed] node [below] {$1$} (m-2-4)
		    edge [dotted, thick] node [right] {$0$} (m-1-3);
	\end{tikzpicture}
	\end{center}
	
	In Section \ref{secsystole}, we discuss the various systole maps and prove the coarse inclusions indicated by solid arrows 	\begin{tikzpicture}
	      \node (1) {};
	      \node (2) [right of=1] {};
	\path[-stealth]
	  (1) edge (2);
      \end{tikzpicture}
       in Corollary \ref{nearsys}. We prove in Section \ref{sechausdist} the inclusions indicated by the double arrows
      \begin{tikzpicture}
	      \node (1) {};
	      \node (2) [right of=1] {};
	\path[-stealth]
	  (1) edge [double] (2);
      \end{tikzpicture}
      in Corollary \ref{sysnearvertex} and Proposition \ref{nearsystole}; and dashed arrows 
      \begin{tikzpicture}
	      \node (1) {};
	      \node (2) [right of=1] {};
	\path[-stealth]
	  (1) edge [dashed] (2);
      \end{tikzpicture}
       in Lemma \ref{nearcyl} using a wide cylinder theorem of Vorobets \cite{Vorobets-cyl}. 
       Our arguments yield explicit constants $\hdSC(g) = O(2^{32g})$.
      The inclusions indicated by dotted arrows
       \begin{tikzpicture}
	      \node (1) {};
	      \node (2) [right of=1] {};
	\path[-stealth]
	  (1) edge [thick, dotted] (2);
      \end{tikzpicture}
      essentially follow by definition.
                           
      In Section \ref{secpolygon}, we describe the construction of the auxiliary polygon $P_q(\alpha)$ and prove some basic properties. These shall be utilised in Section \ref{secbalance} to prove Theorem \ref{thmbalance}. Finally, we discuss some connections between balance points and curve decompositions in Section \ref{secvertexbalance}.
      
      \subsection*{Convention for constants}
      
      We shall label constants as follows:
      \begin{itemize}
       \item $\sK_i$ -- coarse Lipschitz constant,
       \item $\sQ_i$ -- quasiconvexity constant,
       \item $\sP_i$ -- nearest point projection approximation constant,
       \item $\sh_i$ -- bound on Hausdorff distance,
      \end{itemize}
      where $i = 1$ when the constants relate to straight vertex cycles; $i=2$ for systoles; and $i=3$ for cylinders.
      
      \subsection*{Acknowledgements}

	The authors would like to thank Saul Schleimer, Kasra Rafi, and Yair Minsky for interesting conversations. We are greatly appreciative of Jing Tao and Max Forester for providing helpful suggestions during the writing process. We also thank Sam Taylor for encouraging us to find universal bounds on our constants, and Ser-Wei Fu for pointing out Lemma \ref{lengthmin}.
	
	\section{Background}\label{secbackground}
	
	We begin by establishing some definitions and conventions. Throughout this paper, $S$ will be a closed, connected, orientable surface of genus $g \geq 2$.    
	
	\subsection{Coarse geometry}
	
	A constant is called \emph{uniform} if it depends only on the topological type of a surface $S$. A constant is called \emph{universal} if it can be chosen independently of $S$.	
		Given $x,y \in \R$ and constants $\sK, \sC$ write
		\begin{itemize}
		\item $x \prec_{\sK,\sC} y \iff x \leq \sK y + \sC$,
		\item $x \asymp_{\sK,\sC} y \iff x \prec_{\sK,\sC} y$ and $y \prec_{\sK,\sC} x$.
						\end{itemize}
	If $\sK$ and $\sC$ are universal constants, we will also write 
	$x \prec y$ and $x \asymp y$ respectively.
	If $X$ and $Y$ are subsets of a metric space $(\X, d)$, write 
	\begin{itemize}
	\item $X \subseteq_r Y \iff X \subseteq N_r(Y)$,
	\item $X \approx_r Y \iff X \subseteq_r Y$ and $Y \subseteq_r X$.
	\end{itemize}
	Here, $N_r(Y)$ denotes the closed $r$--neighbourhood of $Y$ in $\X$.
	The \emph{Hausdorff distance} between closed sets $X$ and $Y$ in $\X$ is the smallest $r \geq 0$ such that $X \approx_r Y$.
	
	Two relations $f_1 : \X \rightarrow \Y$ and $f_2 : \X \rightarrow \Y$ between metric spaces are said to \emph{agree up to error} $\sh$ if $f_1(x) \approx_\sh f_2(x)$ for all $x \in \X$. We also say that $f_1$ and $f_2$ \emph{coarsely agree}.
	
	Let $I$ be either an interval in $\R$, or the intersection of an interval with $\Z$. A relation $f : I \rightarrow \X$ is a \emph{$(\sK, \sC)$--quasigeodesic} if
	\[\diam_\X(f(s) \cup f(t)) \asymp_{\sK,\sC} |s-t|\]
	for all $s,t \in I$. We say $f : I \rightarrow \X$ is a \emph{reparameterised quasigeodesic} if there is a homeomorphism $h : \R \rightarrow \R$ and constants $\sK, \sC$ and $\sR$ such that $f \circ h$ is a $(\sK, \sC)$--quasigeodesic, and
	\[\diam_\X(f([h(t), h(t+1)])) \leq \sR\]
	for all $t\in h^{-1}(I)$.
	
	A relation $f : \X \rightarrow \Y$ between metric spaces is 	called $(\sK,\sC)$--\emph{coarsely Lipschitz} if for every $x,y\in\X$, 
	\[\diam_\Y(f(x) \cup f(y)) \prec_{\sK,\sC} d_\X(x,y).\]
	If, in addition, there is a constant $\sR \geq 0$ such that $d_\X(y, f(y)) \leq \sR$ for all $y \in \Y$ then we say $f$ is a $(\sK,\sC,\sR)$--\emph{coarse Lipschitz retract} onto $Y$. We shall refer to $(\sK,\sK,\sK)$--coarse Lipschitz retracts as $\sK$--coarse Lipschitz retracts.
	
	A geodesic space $\X$ is \emph{$\hyp$--hyperbolic} if every geodesic triangle is $\hyp$--slim. That is, for all $x,y,z \in \X$ and geodesics $[x,y]$, $[x,z]$ and $[z,y]$, we have
	\[[x,y] \subseteq_\hyp [x,z] \cup [z,y].\]
		
	A subset $\Y$ of a geodesic space $\X$ is \emph{$\sQ$--quasiconvex} if every geodesic in $\X$ connecting a pair of points in $\Y$ is contained in the $\sQ$--neighbourhood of $\Y$.
	
	\begin{Thm}[\cite{Minsky-qc} Lemma 3.3]\label{retractqc}
	 Let $\X$ be a $\hyp$--hyperbolic space and suppose a subset $\Y \subseteq \X$ admits a $(\sK,\sC,\sR)$--coarse Lipschitz retract $f : \X \rightarrow \Y$. Then $\Y$ is $\sQ$--quasiconvex, where $\sQ = \sQ(\hyp,\sK,\sC,\sR)$. \halmos
	\end{Thm}
	
	It is worth noting that a coarse Lipschitz retract to a subset $\Y$ of a $\hyp$--hyperbolic space $\X$ need not coarsely agree with a nearest point projection to $\Y$.
	
	\subsection{Curve graphs}

	By a \emph{curve} on a surface $S$, we shall mean a free homotopy class of simple closed curves which are not null-homotopic. The \emph{curve graph} $\C(S)$ of $S$ is the graph whose vertices are curves on $S$, and whose edges span pairs of distinct curves which have disjoint representatives. We shall also write $\alpha \in \C(S)$ to mean that $\alpha$ is a curve on $S$.
	
	Given curves $\alpha, \beta \in \C(S)$, define their distance $d_S(\alpha,\beta)$ to be the length of a shortest edge-path connecting them in $\C(S)$. The curve graph is locally infinite, and has infinite diameter. Let
	\[i(\alpha,\beta) = \min \{|a \cap b| : a \in \alpha, b \in \beta \}\]
	denote the \emph{geometric intersection number} of $\alpha$ and $\beta$. Combining the distance bounds due to Hempel \cite{Hempel-cc} for $g = 2,3$, and Bowditch \cite{bhb-unif} for $g \geq 4$ gives the following.
	
	\begin{Lem}\label{distbound}
	 Let $\alpha$ and $\beta$ be curves in $\C(S)$ where $i(\alpha,\beta) \geq 1$. Then
	 \[d_S(\alpha,\beta) \prec_{2,2} \log_\Gen i(\alpha,\beta), \]
	 where $\Gen = \max\{2, g-2\}$. \halmos
	\end{Lem}
	
	It follows that $\C(S)$ is connected.
	We shall be using the following observation several times throughout this paper to obtain universal bounds on distances.
	
	\begin{Lem}\label{logmax}
	 For any $C > -2$, the quantity $\log_\Gen(g + C)$ achieves its maximum at $g=4$ among all integers $g \geq 2$. \halmos
	\end{Lem}
	
	 Masur and Minsky prove the following fundamental result concerning the large-scale geometry of curve graphs.
	
	\begin{Thm}[\cite{MM1}]\label{hyperbolic}
	 The curve graph $\C(S)$ is $\hyp$--hyperbolic for some $\hyp > 0$. \halmos
	\end{Thm}
	
	Recent independent work of several authors show that there is a universal hyperbolicity constant \cite{Aougab-unif}, \cite{bhb-unif}, \cite{CRS-unif}, \cite{HPW-unicorn}.

	\subsection{Train tracks}\label{sectrain}
	We now review some background on train tracks.
	For a thorough introduction, see \cite{pennerharer} and \cite{Mosher-train}.
		
	A \emph{pretrack} $\tau$ is a properly embedded graph on $S$ satisfying the following conditions. The edges of $\tau$, called \emph{branches}, are smoothly embedded. Each vertex $s$ of $\tau$, called a \emph{switch}, is equipped with a preferred tangent $v_s \in T_s^1 S$. A half-branch incident to $s$ is called incoming (or outgoing) if its derivative at $s$ is $v_s$ (or $-v_s$). For every switch $s$ we require each incident half-branch to be either incoming or outgoing, and that there is at least once of each type incident to $s$. We allow pretracks to have embedded loops as connected components.
	
	A \emph{train-route} (or \emph{train-loop}) on $\tau$ is a smoothly immersed path (or loop) on $S$ contained in $\tau$. Say a simple closed curve $\alpha\in\C(S)$ is \emph{carried} by $\tau$, and write $\alpha \prec \tau$, if it can be realised as a train-loop on $\tau$. A pretrack $\sigma$ is \emph{carried} by $\tau$, also written $\sigma \prec \tau$, if it can be smoothly homotoped into $\tau$.
	
	Let $\tau$ be a pretrack on $S$. Any curve $\alpha$ realised as a train route on $\tau$ induces a function ${w_\alpha : \B(\tau) \rightarrow \R_{\geq 0}}$ on the branches of $\tau$ defined by setting $w_\alpha(b)$ to be the number of times $\alpha$ runs over the branch $b\in \B(\tau)$. In general, a function $w : \B(\tau) \rightarrow \R_{\geq 0}$ is called a \emph{transverse measure} on $\tau$ if it satisfies the \emph{switch conditions}: at every switch $s$ of $\tau$, we require
	\[\sum_{b \in \B^+(s)} w(b) = \sum_{b \in B^-(s)} w(b), \]
	where $\B^+(s)$ and $\B^-(s)$ are respectively the sets of incoming and outgoing half-branches incident to $s$.
	
	Integral measures on $\tau$ are in one-to-one correspondence with realisations of multicurves on $\tau$, where we allow for multiple parallel copies of each component and trivial loops. One can recover a multicurve from an integral measure as follows: for each branch $b$, take $w(b)$ parallel strands within a \emph{tie neighbourhood} of $b$ on $S$. At each switch $s$, glue the incoming strands to the outgoing strands in pairs -- there is exactly one pairing of the strands which does not lead to self-intersections. See Section 3.7 of \cite{Mosher-train} for more details.
		
	A pretrack $\tau$ on $S$ is called a \emph{train track} if all its complementary regions have negative index: on closed surfaces, this is equivalent to saying that no complementary region is a nullgon, monogon, bigon or annulus. This condition guarantees that every train-loop on a train track $\tau$ is essential on $S$, and that two train-loops on a train track $\tau$ are freely homotopic on $S$ if and only if they agree up to reparameterisation. In particular, distinct integral measures on $\tau$ correspond to distinct homotopy classes of weighted multicurves.
	
	The set $\bM(\tau)$ of transverse measures on a train track $\tau$ defines a convex cone in $\R_{\geq 0}^{\B(\tau)}$ with finitely many extreme rays. 	Each extreme ray of $\bM(\tau)$ contains a unique minimal integral transverse measure on $\tau$, and such a measure arises from a simple closed curve carried by $\tau$.
	We call such curves \emph{vertex cycles} of $\tau$, and write $V(\tau) \subset \C(S)$ for the set of vertex cycles of $\tau$. (Vertex cycles are so named because they form vertices of the polytope obtained by projectivisig $\bM(\tau)$.)
	
	\begin{Thm}[\cite{Mosher-train} Proposition 3.11.3]\label{vertexchar}
	 A simple closed curve $\alpha$ is a vertex cycle of a train track $\tau$ if and only if it realisable as either an embedded loop, a figure--8, or a barbell on $\tau$. \halmos
	\end{Thm}
	
	Every transverse measure $w \in \bM(\tau)$ can be written as a sum of transverse measures arising from vertex cycles of $\tau$. Following a curve surgery argument of Masur, Mosher, and Schleimer with some extra care, one can show the following.
	
	\begin{Thm}\cite{MMS-train}\label{vertexdecomp}
	 Let $\alpha \in \C(S)$ be a multicurve carried by a train track $\tau$. Then there are integers $m_v\geq 0$ for each $v\in V(\tau)$ such that $w_\alpha = \frac{1}{2}\sum_{v \in V(\tau)} m_v w_v$.
	  \halmos
	\end{Thm}

	Aougab shows that $V(\tau)$ has universally bounded diameter in $\C(S)$ by controlling the intersection number between any two curves carried by $\tau$ which run over each branch at most twice.
	
	\begin{Thm}\cite{Aougab-unif}\label{vertexdiam}
	 There is a universal constant $\lipV$ such that for any train track $\tau$ on $S$, the set $V(\tau)$ has diameter at most $\lipV$ in $\C(S)$. \halmos
	\end{Thm}

	If one follows Aougab's proof using Bowditch's distance bound \cite{bhb-unif}, one can show that taking $\lipV = 14$ suffices (even for surfaces with punctures). Masur and Minsky conjecture that $\lipV$ can be taken to be equal to 3, and Aougab shows that this is true for sufficiently large genus.

	A train track is called \emph{generic} if all switches are trivalent. Any train track can be perturbed to a generic one which carries exactly the same curves. Suppose $\tau$ is a generic train track with at least one switch. Every switch of $\tau$ separates its three incident half-branches into one \emph{large} half-branch on one side and a pair of \emph{small} half-branches on the other. A branch $b$ is called \emph{large} (or \emph{small}) if both of its half-branches are large (or small); if it has one of each type then we say it is \emph{mixed}.

	\begin{figure}[h]
	    \centering{
	    \resizebox{100mm}{!}
	    {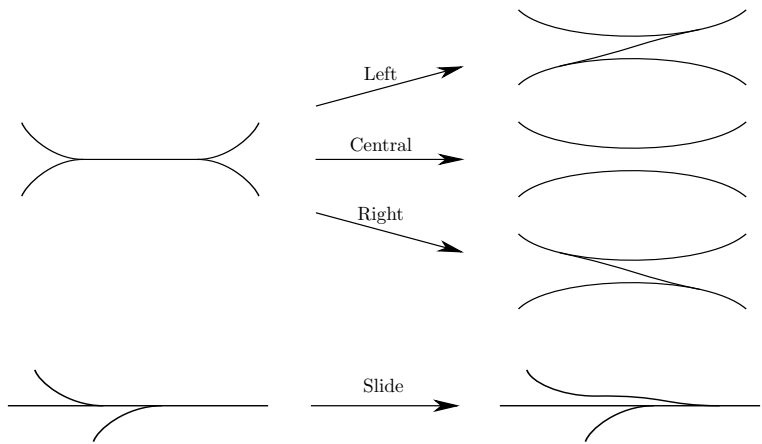}
	    \caption{Splitting a large branch; and sliding a mixed branch.}
	    \label{fig:split}
	    }
	\end{figure}

	We now described splitting and sliding -- refer to Sections 3.11 and 3.12 of \cite{Mosher-train} for more details.

	To a generic train track $\tau$, we may \emph{split} along a large branch or \emph{slide} along a mixed branch to obtain a new train track carried by $\tau$ -- see Figure \ref{fig:split}. We say $(\tau_i)_i$ is a \emph{splitting and sliding sequence of train tracks} if each $\tau_{i+1}$ is obtained from $\tau_i$ by a split or slide. If $\alpha$ is a curve carried by $\tau$, we may \emph{split $\tau$ towards $\alpha$} by splitting a large branch of $\tau$ to obtain a new train track carrying $\alpha$, where the choice of left, right or central split is determined by the transverse measure on $\tau$ induced by $\alpha$. (More specifically, $\alpha$ will either be carried by exactly one of the left or right splits; or both in which case we can choose the central split.) A splitting sequence carried out in this manner will eventually terminate in a train track with an embedded loop component isotopic to $\alpha$.

	\begin{Thm}\cite{MM3}\label{splitting}
	 Let $(\tau_i)_i$ be a splitting and sliding sequence of train tracks. Then $(V(\tau_i))_i$ forms a reparameterised quasigeodesic in $\C(S)$. \halmos
	\end{Thm}

	Hamenst\"adt gives an alternative proof of this result which achieves universal quasiconvexity constants, and also removes the birecurrence assumption for splitting sequences \cite{ham-disc}.

	\subsection{Half-translation surfaces and Teichm\"uller discs}

	A \emph{half-translation} structure $q$ on a closed surface $S$ consists of a finite set of singular points together with an atlas of charts to $\CC \cong \R^2$ defined away from the singular set, where the transition maps are \emph{half-translations}, i.e. of the form $z \mapsto \pm z + c$ where $c \in \CC$. The singular points have Euclidean cone angle of the form $k\pi$ where $k \geq 3$. The atlas determines a preferred vertical direction on $q$. Half-translation surfaces can be constructed by taking a finite collection of disjoint Euclidean polygons with edges glued in pairs isometrically via half-translations.

	The space of half-translations structures on $S$ can be naturally identified with $\QD(S)$, the space of \emph{holomorphic quadratic differentials} on $S$ up to isotopy. A holomorphic quadratic differential on a Riemann surface $X$ is a differential locally of the form $q = q(w)dw^2$, where $q(w)$ is a holomorphic function in the local co-ordinate $w$. The corresponding half-translation structure can be obtained by defining charts near each regular point $z_0$ as 
	\[z \mapsto \int_{z_0}^z \sqrt{q(w)}dw. \]
	An order $k$ zero of $q$ corresponds to a singularity with cone angle $(k+1)\pi$. 	
	A \emph{marking} of a surface $X$ is a homeomorphism $f : S \rightarrow X$ from the reference surface $S$. The \emph{Teichm\"uller space} $\Teich(S)$ of $S$, the space of marked conformal (or complex) structures on $S$ up to isotopy, is homeomorphic to an open ball $\R^{6g - 6}$. The projection map $\QD(S) \rightarrow \Teich(S)$ defined by taking $q\in \QD(S)$ to its underlying conformal structure can be canonically identified with the cotangent bundle to $\Teich(S)$. Restricting this projection to the space $\QD^1(S)$ of \emph{unit-area} half-translation structures on $S$ gives the unit cotangent bundle to $\Teich(S)$. Adopting the convention of Duchin, Leininger, and Rafi \cite{DLR-length}, we will write $q$ to denote both a half-translation surface in $\QD(S)$, as well as its underlying 	conformal structure in $\Teich(S)$. 

	There is a natural $\SL(2,\R)$--action on $\QD(S)$ defined as follows. Let $q\in\QD(S)$ be a half-translation surface and $A \in \SL(2,\R)$ be a real linear transformation on $\R^2$. The half-translation surface $A \cdot q$ has as its atlas the charts obtained by postcomposing each co-ordinate chart of $q$ to $\CC \cong \R^2$ with $A$. One can perform this action by deforming a defining set of polygons for $q$ by $A$ and observing that the gluing patterns are preserved. Also note that $\SL(2,\R)$--deformations preserve area, and so this action descends equivariantly to $\QD^1(S)$. 	
	We shall equip $\Teich(S)$ with the \emph{Teichm\"uller metric} which measures the amount of quasiconformal distortion between two conformal structures. With this metric, $\Teich(S)$ is a complete geodesic space in which geodesics, called \emph{Teichm\"uller geodesics}, are projections of orbits in $\QD(S)$ under the \emph{Teichm\"uller geodesic flow} (or \emph{diagonal action})
		  \[g_t = \begin{pmatrix}
	            e^{t} & 0 \\
	            0 & e^{-t}
	           \end{pmatrix} \in \SL(2,\R)\]
       for $t\in\R$. In particular, $t \mapsto g_t \cdot q$ gives a unit speed paramaterisation of a Teichm\"uller geodesic $\G_q$ through $q$.
	
	A \emph{Teichm\"uller disc} $\disc = \disc_q$ is the projection of the $\SL(2,\R)$--orbit of a half-translation surface ${q \in \QD^1(S)}$ to $\Teich(S)$. Using the Teichm\"uller metric, $\disc$ is an isometrically embedded copy of the hyperbolic plane $\Hy^2 \cong \SO(2,\R) \backslash \SL(2,\R)$ with curvature $-4$ in $\Teich(S)$. We shall write $d_\disc$ for the Teichm\"uller metric restricted to $\disc$. 	

	A \emph{saddle connection} on a half-translation surface is a geodesic segment which meets the singular set precisely at its endpoints.  Let $\alpha \in \C(S)$ be a curve and $q \in \QD(S)$ a half-translation surface. Then either there is a unique maximal flat cylinder foliated by closed geodesic leaves isotopic to $\alpha$, or $\alpha$ has a (unique) geodesic representative on $q$ which is a concatenation of saddle connections. We call $\alpha$ a \emph{cylinder curve} in the former case.

	Write $\alpha^q$ for a geodesic representative of $\alpha$ on $q$, and define the \emph{flat length} of $\alpha$ on $q$ to be $l_q(\alpha) = l(\alpha^q)$. We can also define the \emph{horizontal} and \emph{vertical} lengths of $\alpha$ on $q$, denoted $l_q^H(\alpha)$ and $l_q^V(\alpha)$ respectively, by integrating $\alpha^q$ against the pullbacks of the infinitesimal metrics $|dx|$ and $|dy|$ on $\R^2$.
	Call $\alpha$ \emph{completely horizontal} (or \emph{completely vertical}) on $q$ if $l_q^V(\alpha) = 0$ (or $l_q^H(\alpha) = 0$). We say that a curve $\alpha$ has \emph{constant slope} on $q$ if it is completely horizontal on $\rho_\theta^{-1} \cdot q$ for some rotation $\rho_\theta \in SO(2,\R)$; otherwise it has \emph{non-constant slope} on $q$.

	Observe that
	\[l_{g_t \cdot q}^H(\alpha) = e^tl_q^H(\alpha) \qquad \textrm{and} \qquad l_{g_t \cdot q}^V(\alpha) = e^{-t}l_q^V(\alpha). \]
	Say $\alpha$ is \emph{balanced} on $q$ if $l_q^H(\alpha) = l_q^V(\alpha)$. We call $t \in \R$ the \emph{balance time} of $\alpha$ along a Teichm\"uller geodesic $\G = \G_q$ if $\alpha$ is balanced on $g_t \cdot q$. The balance time exists, and is unique, whenever $\alpha$ is neither completely vertical nor horizontal along $\G$. Define $l_\G(\alpha) = \inf_{q \in \G} l_q(\alpha)$.
	
	\begin{Lem}\label{lengthconvex}
	The flat length of a curve $\alpha$ is strictly convex along any Teichm\"uller geodesic, and satisfies
	\[l_{g_t\cdot q}(\alpha) \leq \sqrt{2}e^tl_q(\alpha) \]
	for all $q \in \QD(S)$.
	 If $\alpha$ is neither completely horizontal nor completely vertical on $q$ then 
	 \[l_{g_t \cdot q}(\alpha) \asymp l_\G(\alpha) \cosh\left(t - t_\alpha\right) \]
	 	 	 	 where $t_\alpha$ is the balance time $\alpha$ along $\G = \G_q$. 	 Moreover, if the (unique) minimum of $l_{g_t \cdot q}(\alpha)$ occurs at $t = t_{m}$ then
	 $|t_{m} - t_\alpha| \leq \cosh^{-1} 2$. \halmos
	\end{Lem}

	\begin{Lem}\label{lengthmin}
	 If $\alpha$ has non-constant slope on $q$, then the flat length of $\alpha$ has a unique minimiser on $\disc_q$.
	\end{Lem}
	
	\proof
	Since $\alpha$ has non-constant slope, its flat length diverges to infinity as we approach the ideal boundary of $\disc_q$. Therefore a minimum must be attained. The minimiser must be unique, for otherwise we can connect distinct minimisers in $\disc_q$ by a Teichm\"uller geodesic, a contradiction.
	\endproof
	
		For a Teichm\"uller disc $\disc$, the boundary circle $\partial \disc \subset \PMF(S)$ corresponds to the measured foliations on $S$ which can be realised as a constant slope foliation on some $q \in \disc$, where the transverse measure is the Lebesgue measure (up to scale).

	\section{Train tracks and multi-arcs on flat surfaces}\label{secstraight}
	
	  \subsection{Train tracks induced by geodesics on flat surfaces}
 
	  Throughout this section, fix a half-translation surface ${q \in \QD^1(S)}$. Suppose a curve $\alpha\in\C(S)$ is not a cylinder curve on $q$, so that it has a unique geodesic representative $\alpha^q$. We construct a train track $\tq(\alpha)$ on $S$ by performing a ``smoothing'' operation at each singular point. 
	  	  
	 \begin{figure}[h]
	    \centering{
	    	    {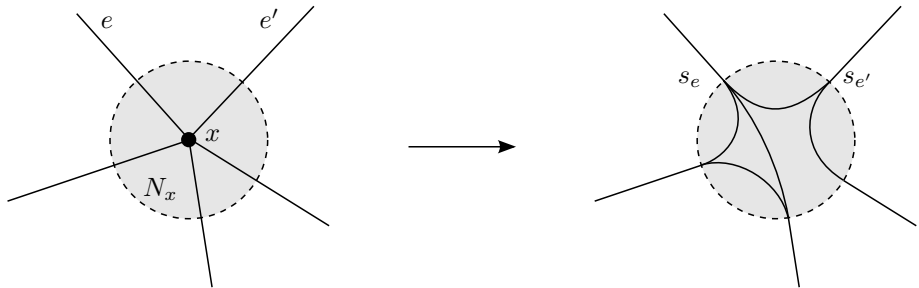}
	    \caption{Smoothing $\alpha^q$ at a singularity $x$ to locally produce a train track.}
	    \label{fig:smooth}
	    }
	\end{figure}
  
	  Let $\Gamma = \Gamma_q(\alpha)$ be the embedded graph on $q$ whose vertices and edges are respectively the singular points and saddle connections used by $\alpha^q$. At each singular point $x$, delete from $\Gamma$ a small open regular neighbourhood $N_x$. The \emph{link} $\Gamma \cap \partial N_x$ of $x$ is in natural bijection with the set of half-edges incident to $x$ -- write $s_e \in \Gamma \cap \partial N_x$ for the point associated to a half-edge $e$.
	  Add a smoothly embedded arc, with interior inside $N_x$, connecting two such points $s_e$ and $s_{e'}$ if and only if $\alpha^q$, regarded as a closed path on $q$, passes through $x$ by entering along $e$ and exiting along $e'$ (or vice versa) -- see Figure \ref{fig:smooth}. (Note that $e$ and $e'$ are necessarily distinct as $\alpha^q$ is a geodesic.) Such an arc can be realised so that at each endpoint, the incoming unit tangent vector along the arc coincides with the outgoing unit tangent vector along the corresponding half-edge. This gives a switch structure to the points $s_e$ which are the endpoints of at least two new arcs.
	  Since $\alpha$ is simple, we can arrange so that all the new arcs have disjoint interiors.

	  The above procedure produces a pretrack $\tq(\alpha)$ whose switches are the points $s_e$ with at least three incident arcs. It may be that $\tq(\alpha)$ has no switches and is thus an embedded loop representing $\alpha$. If $\alpha$ is a cylinder curve on $q$, we take $\tq(\alpha)$ to be any closed geodesic leaf representing $\alpha$ not passing through any singular points. We can also naturally extend the definition of $\tq(\alpha)$ to multicurves, noting that it could possibly be disconnected.
	  
	  A multicurve $\alpha$ can be realised as a train-loop on $\tq(\alpha)$ by perturbing it near the singularities: each time the closed path $\alpha^q$ passes through a singularity $x$ along $e$ and $e'$, perturb it within $N_x$ so that it instead runs along the arc on $\tq(\alpha)$ connecting $s_e$ and $s_{e'}$.
	  
	  Conversely, one can ``straighten'' $\tq(\alpha)$ to $\Gamma$: whenever there is an arc of $\tq(\alpha)$ connecting $s_e$ and $s_{e'}$ within some $N_x$, homotope it to a geodesic path running along $e$ and $e'$ via $x$. This homotopy sends any train-route $\eta$ on $\tq(\alpha)$ to the unique geodesic path on $q$ connecting the endpoints of $\eta$ in its relative homotopy class. In particular, any train-loop on $\tq(\alpha)$ is sent to a geodesic representative which lies completely within $\Gamma$. Note that $\tq(\alpha)$ is identical to its straightening if and only if every component of $\alpha$ is a cylinder curve.

	  \begin{Lem}\label{tcarry}
	   Let $\alpha\in\C(S)$ be a multicurve. Then $\tq(\alpha)$ is a train track carrying $\alpha$. Furthermore, if $\beta$ is a multicurve carried by $\tq(\alpha)$ then $\tq(\beta)$ is a subtrack of $\tq(\alpha)$. 	  \end{Lem}
	
	  \proof
	  To show that $\tq(\alpha)$ is a train track we need to verify that complementary region of $\tq(\alpha)$ on $S$ has negative index. The remaining claims follow from the construction.
  
	  Suppose $Y$ is a complementary region of $\tq(\alpha)$ of one of the following forbidden types. We apply the above ``straightening'' procedure to $\partial Y$. 	 	  If $Y$ is a nullgon or monogon, then $\partial Y$ is a closed train-route which straightens to a closed geodesic on $q$ homotopic to a constant path. 	  	  If $Y$ is a bigon, then $\partial Y$ consists of two distinct train-routes connecting a pair of switches $s$ and $s'$. Upon straightening, these train routes become distinct geodesics connecting $s$ and $s'$ in the same relative homotopy class. Either of these scenarios contradict the fact that there is a unique geodesic connecting a given pair of points on $q$ in each relative homotopy class. 	  
	  Now suppose $Y$ is an annulus. Then $\partial Y$ consists of a pair of homotopic train-loops, say $\gamma$ and $\gamma'$, which straighten to a pair of distinct homotopic closed geodesics. Since these closed geodesics are distinct, they must be core curves of a common flat cylinder on $q$. The integral measure $w_\alpha$ on $\tq(\alpha)$ determined by $\alpha$ places positive weight on each branch -- in particular, all branches which $\gamma$ and $\gamma'$ run over. As one reconstructs $\alpha$ from $w_\alpha$ using the process described in Section \ref{sectrain}, we see that $\alpha$ contains two components which are parallel to $\gamma$ and $\gamma'$ respectively. Since $\gamma$ and $\gamma'$ are homotopic, this gives a contradiction.
	  \endproof

	 \subsection{Straight vertex cycles for Teichm\"uller discs}
	  
	  The construction of $\tq(\alpha)$ is completely determined by the sequence of directed saddle connections used by $\alpha^q$. Since this data is preserved under $\SL(2,\R)$--deformations on $q$, we have the following observation.
	  
	  \begin{Lem}
	   For any half-translation surface $q \in \QD^1(S)$, multicurve $\alpha \in \C(S)$ and $A \in \SL(2,\R)$, the train track $A \cdot \tq(\alpha)$ is isotopic to $\tau_{A \cdot q}(\alpha)$ on $A \cdot q$. \halmos
	  \end{Lem}
	  
	  This allows us to extend the construction to Teichm\"uller discs. Let $\disc \subset \Teich(S)$ be a Teichm\"uller disc and $\alpha \in \C(S)$ a multicurve. Choose any $q \in \disc$, and let $f_q : S \rightarrow q$ be a marking. Then
	  \[\td(\alpha) := f_q^{-1}(\tq(\alpha))\]
	  is a train track on $S$. By the above lemma, the isotopy class of $\td(\alpha)$ on $S$ is independent of the choice of $q\in\disc$.
	  In addition, Lemma \ref{tcarry} also holds if one replaces $\tq$ with $\td$.
	  
	  Given a Teichm\"uller disc $\disc \subset \Teich(S)$, define its set of \emph{straight train tracks}
	  \[\TT(\disc) := \{\td(\alpha) ~|~ \alpha \textrm{ a multicurve on } S \} \subset \TT(S)\]
	  and its \emph{straight vertex cycle} set
	  \[V(\disc) := \bigcup_{\tau \in \TT(\disc)} V(\tau). \]
	  Appealing to the characterisation of vertex cycles, straight vertex cycles for $\disc$ have geodesic representatives on any $q \in \disc$ which run over each saddle connection at most once in each direction. One can then think of $V(\disc)$ as the set of combinatorially short curves with respect to $\disc$.	   
	  We shall also write $\VD(\alpha) = V(\td(\alpha))$.
	  
	  Let $\lipV$ be a universal upper bound on the diameter of the vertex set of any train track.
	  
	  \begin{Lem}\label{vertexmap}
	   The relation $\VD : \C(S) \rightarrow V(\disc)$ is a $\lipV$--coarsely Lipschitz retract. 	  \end{Lem}
	  
	  \proof
	  Let $\alpha$ and $\beta$ be disjoint curves. By the Lemma \ref{tcarry}, $\td(\alpha)$ and $\td(\beta)$ are both subtracks of $\td(\alpha \cup \beta)$, from which it follows that $\VD(\alpha) \cup \VD(\beta) \subseteq \VD(\alpha \cup \beta)$. Since the diameter of $\VD(\alpha \cup \beta)$ in $\C(S)$ is at most $\lipV$, it follows that $\VD$ is {$\lipV$--coarsely} Lipschitz.
	  
	  Now suppose $\gamma \in V(\disc)$, which means there is some multicurve $\alpha$ such that $\gamma \in V(\td(\alpha))$. By Lemma \ref{tcarry}, $\td(\gamma)$ is a subtrack of $\td(\alpha)$ which also carries $\gamma$. By the characterisation of vertex cycles, $\gamma$ is also a vertex cycle of $\td(\gamma)$. Therefore $\gamma \in \VD(\gamma)$ and so $\VD$ is a coarse retract.
	  \endproof
	  
	  It immediately follows, via Theorem \ref{retractqc}, that $V(\disc)$ is quasiconvex.
	  
	  \begin{Cor}
	   There is a universal constant $\qcV$ such that for any Teichm\"uller disc $\disc \subset \Teich(S)$, the set $V(\disc)$ is $\qcV$--quasiconvex in $\C(S)$. \halmos
	  \end{Cor}
	  
	  \begin{Prop}\label{vertexproj}
	   The operation $\VD : \C(S) \rightarrow V(\disc)$ agrees with the nearest point projection map from $\C(S)$ to $V(\disc)$ up to a universal error $\projV$. 	  \end{Prop}
	  
	  \proof
	  Perturb $\td(\alpha)$ to a generic train track $\tau_0$, and consider a train track splitting and sliding sequence 
	  $\tau_0 \succ \tau_1 \succ \ldots \succ \tau_n = \alpha$.
	  Choose a vertex cycle $\gamma_i \in V(\tau_i)$ for $0 \leq i \leq n$.
	  By Theorem \ref{splitting}, the sequence $(\gamma_i)_i$ forms a reparameterised quasigeodesic from $\VD(\tau)$ to $\alpha$ in $\C(S)$ with universal quasiconvexity constants. 
	  
	  Let $\beta$ be a nearest point projection of $\alpha$ to $V(\disc)$. Since $V(\disc)$ is $\qcV$--quasiconvex, $\beta$ must lie within a universal distance of any geodesic connecting $\alpha$ to $\VD(\alpha)$, and hence the quasigeodesic $(\gamma_i)_i$. Therefore, there is some $\gamma_k$ which is universally close to $\beta$ and hence $V(\disc)$. By Lemma \ref{vertexmap}, $\VD$ is a $\lipV$--coarse Lipschitz retract which means $\VD(\gamma_k)$ is universally close to $\gamma_k$. But by Lemma \ref{tcarry}, $\td(\gamma_k)$ is a subtrack of $\td(\alpha)$, and so $\VD(\gamma_k) \subseteq \VD(\alpha)$. It follows that $\beta$ is within a universally bounded distance of $\td(\alpha)$.
	  \endproof
	  
	  \subsection{Filling multi-arc graphs for Teichm\"uller discs}

	 Given $q\in \disc$, we consider the underlying surface $S$ with marked points at the singularities of $q$.
	 More precisely, let $S' = (S,Z)$ be the reference surface $S$ whose set of marked points $Z$ is the preimage of the singular points on $q$ via the marking $f_q : S \rightarrow q$. After isotopy, we may also assume that for every $q' \in \disc$, the marking $f_{q'} : S \rightarrow q'$ maps $Z$ bijectively to the singular points of $q'$. Thus, we may consider the reference surface with marked points $S' = S'(\disc)$ associated to $\disc$.
	 
	 The \emph{arc complex} $\A(S')$ of $S'$ has as vertices proper isotopy classes of embedded arcs on $S$ whose endpoints lie in $Z$. The simplices of $\A(S')$ correspond to \emph{multi-arcs}: a set of arcs which can be properly isotoped to have pairwise disjoint interiors. Let $\A(\disc)$ be the subcomplex of $\A(S')$ spanned by the arcs which can be realised as a saddle connection on some (hence any) $q \in \disc$. The arc complex is connected, has infinite diameter, and is locally infinite. Moreover, its $1$--skeleton is universally hyperbolic \cite{HPW-unicorn}.

	 	 The $1$--skeleton of the first barycentric subdivision of $\A(S')$ is the \emph{multi-arc graph} $\MA(S')$. This naturally contains the subgraph $\MA(\disc)$ spanned by multi-arcs realisable as saddle connections on some $q\in\disc$. Observe that the geodesic representative of any multi-arc (or multicurve) on $q$ forms a collection of disjoint saddle connections, which in turn gives a vertex in $\MA(\disc)$. So there is a $1$--Lipschitz map $\MA(S') \rightarrow \MA(\disc)$ that coincides with the identity when restricted to $\MA(\disc)$.
	 	 
	  \begin{Lem}
	  There is a $1$-Lipschitz retract $\MA(S')\rightarrow \MA(\disc)$. In particular $\MA(\disc)$ is connected and isometrically embedded in $\MA(S')$. \halmos
	  \end{Lem}

	  Let $\arcd : \C(S) \to \MA(\disc)$ be the $1$--Lipschitz map that sends a curve $\alpha$ to the set of saddle connections appearing on $\alpha^q$ for any $q \in \disc$. If $\alpha$ is a cylinder curve, we set $\arcd(\alpha)$ to be the saddle connections appearing on the boundary of the maximal flat cylinder on $q$ with core curve $\alpha$.
	  
	  \begin{Lem}\label{arcfill}
	  Suppose that $\arcd(\alpha)$ does not fill $S'$. Then there exists $\beta\in V(\Delta)$ such that $d_S(\alpha,\beta)\leq 2$.
	  \end{Lem}

	  \proof
	  Choose any $\beta\in \VD(\alpha)$. Then by Lemma \ref{tcarry}, $\td(\beta)$ is a subtrack of $\td(\alpha)$, and so upon straightening it follows that $\arcd(\beta) \subseteq \arcd(\alpha)$ (viewed as sets of saddle connections). Since $\arcd(\alpha)$ does not fill $S'$, there is some $\gamma \in \C(S)$ disjoint from $\arcd(\alpha)$ and hence $\arcd(\beta)$. It follows that $\gamma$ is disjoint from both $\alpha$ and $\beta$ .
	  \endproof

	  Define the \emph{filling multi-arc graph} $\FMA(S')$ to be the subgraph of $\MA(S')$ spanned by filling multi-arcs on $S'$. In contrast with $\MA(S')$, this graph is locally finite. Let $\FMA(\disc)$ be the subgraph of $\FMA(S')$ spanned by the filling multi-arcs realisable as saddle connections on any $q \in \disc$.
	  	  	  Restricting the retraction $\MA(S')\rightarrow \MA(\disc)$ to $\FMA(S')$, we deduce the following.

	  \begin{Lem}\label{fmalip}
	  There is a $1$--Lipschitz retract $\FMA(S')\rightarrow \FMA(\disc)$. In particular, $\FMA(\disc)$ is connected, locally finite, and isometrically embedded in $\FMA(S')$. \halmos
	  \end{Lem}

	  We are now ready to prove a bounded geodesic image theorem for the projection $\arcd$.

	  \begin{Thm}\label{bgit}
	   There are universal constants $\sA$ and $\sB$ so that if $G$ is a geodesic in $\C(S)$ disjoint from the $\sA$--neighbourhood of $V(\disc)$, then $\arcd(G)= \{ \arcd(\alpha) ~|~ \alpha\in G\}$ has diameter at most $\sB$ in $\FMA(\disc)$.
	  \end{Thm}
	  
	  \proof
	  Let $\beta\in G$ be a curve which is closest to $V(\disc)$ among all curves on $G$. We claim that $\arcd(G)$ lies in a universally bounded neighbourhood of $\arcd(\beta)$ in $\FMA(\disc)$.
	  
	  Choose any curve $\alpha$ on $G$, and let $G'$ be a geodesic connecting $\alpha$ to $\VD(\alpha)$ in $\C(S)$.	  
	  Appealing to $\hyp$--hyperbolicity of $\C(S)$ and $\qcV$--quasiconvexity of $V(\disc)$, we deduce that $G'$ must pass within a distance $\sC = \sC(\hyp, \qcV)$ of $\beta$, so long as $G$ does not come within a distance $\sA = \sA(\hyp, \qcV)$ of $V(\disc)$.
	  
	  Now consider the train track $\td(\alpha)$. After perturbing to a generic train track $\tau_0$, we may produce a splitting and sliding sequence $\tau_0 \succ \tau_1 \succ \ldots \succ \tau_n = \alpha$. By Theorem \ref{splitting}, the vertex cycles along this splitting sequence form a quasigeodesic which agrees with $G'$ up to universal Hausdorff distance. It follows that there is some vertex cycle $\gamma$ along this splitting sequence within a universal distance $\sC'$ of $\beta$. Let $G''$ be a geodesic connecting $\gamma$ to $\beta$. Taking $\sA$ to be larger if necessary, we may assume every vertex along $G''$ is at least distance 3 from $V(\disc)$. Thus, by Lemma \ref{arcfill}, every curve on $G''$ maps to a filling multi-arc under $\arcd$. Since $\arcd$ is $1$--Lipschitz, $\arcd(G)$ is a path of length at most $\sC'$ in $\FMA(\disc)$ connecting $\arcd(\gamma)$ to $\arcd(\beta)$.
	  
	  Since $\gamma$ is carried by $\td(\alpha)$, we may apply Lemma \ref{tcarry} and argue as in Lemma \ref{arcfill} to show that $\arcd(\gamma) \subseteq \arcd(\alpha)$.
	  In particular, $\arcd(\gamma)$ and $\arcd(\alpha)$ are equal or adjacent in $\FMA(\disc)$. Thus $d_{\FMA(\disc)}(\arcd(\alpha), \arcd(\beta)) \leq \sC' + 1$, and so $\arcd(G)$ has diameter at most $\sB = 2\sC' + 2$ in $\FMA(\disc)$.
	  \endproof

	  \begin{Rem}
	  This theorem is inspired by, but does not seem to follow readily from, Masur and Minsky's bounded geodesic image theorem for subsurface projections -- see Theorem 3.1 of \cite{MM2}. 
	  \end{Rem}
	  
	  The valency of each vertex of $\mathcal{FMA}(\disc)$ can be bounded in terms of the topology of $S'$. In particular, Theorem \ref{bgit} shows that there is a uniform bound on the number of distinct saddle connections that appear on the geodesic representatives of curves along $G$. Moreover, one can uniformly bound the pairwise intersection number between arcs in $\arcd(G)$. In particular, if the geodesic representatives of curves $\alpha$ and $\beta$ use saddle connections $a \in \arcd(\alpha)$ and $b \in \arcd(\beta)$ such that $i(a,b)$ is sufficiently large, then any geodesic $G$ connecting $\alpha$ and $\beta$ in $\C(S)$ must come close to $V(\disc)$.
	  
	  \begin{Cor}
	   There exists a uniform constant $\sN = \sN(S)$ such that for any geodesic $G$ in $\C(S)$ which does not meet the $\sA$--neighbourhood of $V(\disc)$, we have $\#\{a \in \arcd(\alpha) ~|~ \alpha \in G \} \leq \sN$. \halmos
	   	  \end{Cor}

	  We note that $\arcd(\alpha)$ and $\td(\alpha)$ may also be viewed as \textit{markings} on $S$. (This is not the same notion as a marking from a reference surface $S$.) We may define a marking to be (the isotopy class of) an embedded graph on $S$ that fills, with a fixed uniform bound on the numbers of vertices and edges. It then follows that there are only finitely many such homeomorphism classes of markings on $S$. The \textit{marking graph} $\M(S)$ can be defined by taking the set of markings as vertices, with an edge spanning two markings if they have bounded intersection number (the bound should be taken large enough to ensure that the marking graph is connected). The mapping class group $\mcg(S)$ acts properly and cocompactly on $\M(S)$ by isometries and so, by the \v{S}varc--Milnor Lemma, it is quasi-isometric to the marking graph. The same proof as the above theorem gives:

	  \begin{Cor}
	  Let $G$ be a geodesic in $\C(S)$ disjoint from the $\sA$--neighbourhood of $V(\disc)$. Then the sets $\arcd(G)$ and $\td(G) = \{ \td(\alpha) ~|~ \alpha\in G\}$ have uniformly bounded diameter in $\M(S)$. \halmos
	  \end{Cor}

	\section{Systole sets}\label{secsystole}
	
	A \emph{systole} on a metric surface is an essential curve of shortest length. Given a Teichm\"uller disc ${\disc \subset \Teich(S)}$, define the \emph{flat systole map} to be the relation ${\sys : \disc \rightarrow \C(S)}$ which assigns each flat surface $q \in \Delta$ to its set of systoles. We can similarly define the systole maps $\sys^\Ext$ and $\sys^\Hyp$ on $\disc$ with respect to extremal length and hyperbolic length respectively.
	
	It is well-known that the hyperbolic and extremal systole maps defined on all of $\Teich(S)$ coarsely agree. We show that the flat systole map also coarsely agrees with these maps when restricted to $\disc$.
	
	\begin{Prop}\label{nearsys}
	 There is a universal constant $\hdsys$ such that the maps $\sys$, $\sys^\Ext$ and $\sys^\Hyp$ agree up to error $\hdsys$.
	 Consequently, we have
	 \[\sys(\disc) \approx_\hdsys \sys^\Ext(\disc) \approx_\hdsys \sys^\Hyp(\disc) \approx_\hdsys \sys(\disc),\]
	 where the Hausdorff distance is measured in $\C(S)$.
	 Choosing $\hdsys = 37$ suffices.
 	\end{Prop}
	
	Let us first consider flat systoles. One can use a standard argument by bounding the injectivity radius to show the following.
	
	\begin{Lem}\label{syslength}
	 Any systole on a unit-area half-translation surface has length at most $\frac{2}{\sqrt{\pi}}$. \halmos
	\end{Lem}
	
	The main ingredient in showing that the relation $\sys : \disc \rightarrow \C(S)$ is coarsely well-defined is the existence of annuli of definite width -- a key step in both Masur and Minsky's \cite{MM1}, and Bowditch's \cite{bhb-int} proofs of hyperbolicity of the curve graph. Masur and Minsky used a limiting argument, whereas Bowditch's method can produce explicit bounds assuming only an isoperimetric inequality on a given singular Riemannian surface.
	
	\begin{Prop}[\cite{bhb-unif}]\label{wideannulus}
	 For any $q \in \QD^1(S)$, there is a topological annulus on $q$ whose width is at least $\widthB = \frac{\sqrt{2\pi}}{4(2g-1)(2g+6)} \asymp \frac{1}{g^2}$. \halmos
	\end{Prop}
		
	Following their arguments with extra care in computing the constants, we deduce the following. Recall that $\Gen = \max\{2, g-2\}$.	
	\begin{Lem}\label{shortbounded}
	For any $q\in\QD^1(S)$ and any pair of intersecting curves $\alpha, \beta \in \C(S)$, we have
	\[d_S(\alpha,\beta) \prec_{2, \diamL} \log_\Gen l_q(\alpha) + \log_\Gen l_q(\beta)\]
		for some universal constant $\diamL = 33.2$. In particular, for any $L\geq 0$, the set of curves whose length on $q$ is at most $L$ has diameter in $\C(S)$ at most $4\log_\Gen L +\diamL$.
	\end{Lem}
	
	\proof 
	Let $\gamma$ be the core curve of an annulus of width at least $\widthB$ on $q$. Applying Lemma \ref{distbound}, we get
	\begin{eqnarray*}
	d_S(\alpha,\gamma) &\leq& 2 \log_\Gen i(\alpha,\gamma) + 2\\
	&\leq& 2 \log_\Gen \left(\frac{l_q(\alpha)}{\widthB}\right) + 2\\
	&=& 2 \log_\Gen l_q(\alpha) + 2 \log_\Gen\left(\frac{4(2g-1)(2g+6)}{\sqrt{2\pi}}\right) + 2.
	\end{eqnarray*}
		By Lemma \ref{logmax}, $\log_\Gen\left\{(2g-1)(2g+6)\right\}$ is maximised when $g=4$ among all integers $g \geq 2$ and so
	\begin{eqnarray*}
	d_S(\alpha,\gamma) &\leq& 2 \log_\Gen l_q(\alpha) + 2 \log_2\left(\frac{392}{\sqrt{2\pi}}\right) + 2 \\
	&\leq&  2 \log_\Gen l_q(\alpha) + 16.6.
	\end{eqnarray*}
	Applying the triangle inequality to $d_S(\alpha,\beta)$ completes the proof. 	\endproof

	\begin{Lem}\label{syslip}
	 The systole map $\sys : \disc \rightarrow \C(S)$ is $(\lipS,\lipcS)$--coarsely Lipschitz, where $\lipS = \frac{2}{\log\Gen}$ and $\lipcS = 36$.
	\end{Lem}
	
	\proof
	Given $q,q'\in\disc$, let $\alpha\in\sys(q)$ and $\beta\in\sys(q')$ be respective systoles. After rotation, we may assume $q = g_t \cdot q'$, where $t = d_\disc(q,q')$. By Lemmas \ref{lengthconvex} and \ref{syslength}, we have
			\[l_q(\beta) \leq \sqrt{2}e^tl_{q'}(\beta) \leq \frac{2\sqrt{2}}{\sqrt{\pi}} e^t.\]
	Applying Lemma \ref{shortbounded} gives
	\begin{eqnarray*}
	d_S(\alpha,\beta) &\leq& 2\log_\Gen \left(\frac{2\sqrt{2}}{\sqrt{\pi}}\right) + 2\log_\Gen \left(\frac{2\sqrt{2}}{\sqrt{\pi}} e^t\right) + 33.2\\
	  &\leq& \frac{2t}{\log \Gen} + 2\log_2 \left(\frac{8}{\pi}\right) + 33.2\\
	  &\leq& \frac{2}{\log \Gen}d_\disc(q,q') + 36
	\end{eqnarray*}
	and we are done.
	\endproof
	
	Let $X \in \Teich(S)$ be a conformal class of metrics on $S$. The \emph{extremal length} of $\alpha\in\C(S)$ on $X$ is
	\[\Ext_X(\alpha) = \sup_{\rho \in X} \frac{l_\rho(\alpha)^2}{\area(\rho)}\]
	where $l_\rho(\alpha)$ denotes the geodesic length of $\alpha$ on $\rho$.	
	It immediately follows that for any $q\in\QD^1(S)$, we have $l_q(\alpha) \leq \sqrt{\Ext_X(\alpha)}$, where $X \in \Teich(S)$ is the conformal class of $q$.		
	Let $\Hyp_X(\alpha)$ denote the geodesic length of $\alpha$ on the unique hyperbolic metric in the conformal class $X$. A theorem of Maskit \cite{Maskit-comparison} states that 
	\[\frac{\Hyp_X(\alpha)}{\pi} \leq \Ext_X(\alpha) \leq \frac{1}{2}\Hyp_X(\alpha) e^{\frac{1}{2}\Hyp_X(\alpha)}.\]
	Using standard hyperbolic geometry, any hyperbolic systole $\alpha$ on $X$ satisfies $\Hyp_X(\alpha) \leq 2 \log(4g-2)$.
	
	\begin{Lem}
	If $\alpha$ is an extremal or hyperbolic systole on $X \in \Teich(S)$ then
	\[l_q(\alpha)^2 \leq (4g-2)\log(4g-2)\]
	for any $q \in \QD^1(S)$ in the conformal class $X$. 
	\halmos
	\end{Lem}

	Proposition \ref{nearsys} is a consequence of the following.
	
	\begin{Lem}\label{sysagree}
	Let $q \in \QD^1(S)$. Then
	\[\diam_{\C(S)}\{\sys(q) \cup \sys^\Ext(q) \cup \sys^\Hyp(q) \} \leq 37. \]
	Consequently, the maps $\sys$, $\sys^\Ext$ and $\sys^\Hyp$ are coarsely well-defined and agree up to universally bounded error.	 
	\end{Lem}
	
	\proof
	Fix a surface $q\in\disc$, and choose curves $\alpha,\beta \in \sys(q) \cup \sys^\Ext(q) \cup \sys^\Hyp(q)$.
	By the above corollary and Lemma \ref{shortbounded}, we have
	\[d_S(\alpha,\beta) \leq 2\log_\Gen((4g-2)\log(4g-2)) + \diamL.\]
	A computation shows that this quantity is maximised at $g=4$ among all integers $g \geq 2$, and so
	\[d_S(\alpha,\beta) \leq 2\log_2(14\log14) + \diamL \leq 37.3\]
	which completes the proof.
	\endproof

	\section{Bounding Hausdorff distance}\label{sechausdist}
	
	In this section, we prove that the systole set and the set of vertex cycles for $\disc$ agree up to universally bounded Hausdorff distance through a sequence of coarse inclusions. We will also show that the cylinder set for $\disc$ agrees with the systole set up to uniformly bounded Hausdorff distance.
	
	\subsection{Systoles are near straight vertex cycles}
	
	Fix a Teichm\"uller disc $\disc$ and choose some $q \in \disc$. Any curve $\alpha\in\C(S)$ induces a minimal integral transverse measure $w_\alpha$ on $\tau = \tq(\alpha)$. By Theorem \ref{vertexdecomp}, we can write $w_\alpha = \frac{1}{2}\sum_{v \in V(\tau)} m_v w_v$, for some non-negative integers $m_v$. By straightening $\tq(\alpha)$ and counting the number of times $\alpha^q$ and each $v^q$ run over each saddle connection, we can deduce 
	\[l_q(\alpha) = \frac{1}{2}\sum_{v \in V(\tau)} m_v l_q(v).\]
		
	\begin{Lem}
	 Suppose $\alpha\in\sys(\disc)$. Then either $\alpha$ is a vertex cycle of $\tau = \td(\alpha)$, or there are distinct vertex cycles $v, v' \in V(\tau)$ such that $w_\alpha = \frac{1}{2}(w_v + w_{v'})$. 	\end{Lem}
	
	\proof
	If $\alpha$ is a systole on some $q \in \disc$, then
	\[l_q(\alpha) = \frac{1}{2}\sum_{v \in V(\tau)} m_v l_q(v) \geq \frac{1}{2}\sum_{v \in V(\tau)} m_v l_q(\alpha),\]
	and hence the $m_v$ must sum to at most 2. Since $w_\alpha$ and the $w_v$ are minimal integral transverse measures on $\tau$, we deduce that $\sum_{v \in V(\tau)} m_v = 2$. The desired result then follows.
	\endproof
	
	It follows that any curve $\alpha\in\sys(\disc)$ must run over each branch of $\td(\alpha)$ at most twice. By Theorem \ref{vertexdiam} (see also the surrounding remarks), $\alpha$ is at most a distance $\lipV \leq 14$ from any vertex cycle of $\tau$.
	
	\begin{Cor}\label{sysnearvertex}
	 For any Teichm\"uller disc $\disc$, we have $\sys(\disc) \subseteq_{14} V(\disc)$. \halmos
	\end{Cor}

	\subsection{Straight vertex cycles are near systoles}
	
	Before we prove the reverse inclusion, we first state some technical lemmas.
	
	\begin{Lem}\label{hnhd}
	If $b$ is a vertical arc with no singular points on its interior on a unit-area flat surface $q$, then there exists a horizontal arc of length at most $\frac{1}{l_q(b)}$ whose endpoints lie on $b$ with no singular points on its interior. 
	\end{Lem}

	\proof
	We use a Poincar\'e recurrence argmument for the geodesic flow on $q$ in the horizontal direction.
	
	Define the \emph{horizontal neighbourhood} $N_r^H(b)$ of $b$ with radius $r \geq 0$ to be the set of points $x$ on $q$ for which there exists a horizontal arc of length at most $r$ starting from a point on $b$ and ending at $x$ with no singular points on its interior. Observe that $N_r^H(b)$ is the image of a locally isometric immersion $\iota: R \rightarrow q$, where $R$ is a Euclidean rectangle $[-r, r] \times b$ with finitely many horizontal intervals of the form $(t, r] \times \{p\}$ or $[-r, t) \times \{p\}$ removed. These removed intervals correspond to when a horizontal arc emanating from $b$ hits a singularity and cannot be uniquely extended. Moreover, the immersion is locally area preserving. If $\iota$ is an embedding, then \[\area(N_r^H(b)) = \area(R) = 2rl_q(b).\]
	Since $q$ is assumed to have unit area, $\iota$ cannot be an embedding when $r \geq \frac{1}{2l_q(b)}$. Thus there is a non-singular point $x \in q$ which can be connected to $b$ using two distinct horizontal arcs of length at most $\frac{1}{2l_q(b)}$. Taking the union of the two horizontal arcs produces the desired horizontal arc with endpoints on $b$.
	\endproof

	We introduce \emph{$\disc$--bicorns} as an intermediate step. By definition, these are curves which have representatives formed by taking the union of two straight line segments on some (hence all) $q \in \disc$. Such curves are non-trivial since geodesics segments on a non-positively curved surface are unique in their homotopy class relative to their endpoints.
	
	\begin{Lem}\label{nearbicorn}
	 Let $\Gamma$ be an embedded graph on some $q \in \disc$, with vertices at singular points and whose edges are saddle connections. Then there is a cylinder curve or a $\disc$--bicorn which intersects each edge of $\Gamma$ at most once.
	\end{Lem}
	
	\proof
	Let $e$ be an edge of $\Gamma$, which we may assume to be vertical. Applying the Poincar\'e recurrence argument as in the above lemma, there is a horizontal arc $a$ with endpoints on $e$ with no singularities in its interior. If $a$ intersects each edge of $\Gamma$ at most once, we may concatenate it with a (possibly degenerate) subarc $e'\subseteq e$ to form a $\disc$--bicorn intersecting each edge of $\Gamma$ at most once as desired. If not, we may choose an innermost subarc $a' \subseteq a$ with the property that its endpoints are on the same edge of $\Gamma$. The arc $a'$ intersects each edge of $\Gamma$ at most once, and so we may argue as above.
	\endproof
	
	\begin{Lem}\label{nearshort}
	 Given any $\disc$--bicorn $\beta$, there is some $\gamma\in\C(S)$ with $l_\disc(\gamma) \leq 2$ such that $d_S(\beta,\gamma) \leq 2$.
	\end{Lem}

	\proof
	Applying a suitable $SL(2,\R)$--deformation, we may assume that $\beta$ has a bicorn representative on some $q\in\disc$ where the two line segments are in the horizontal and vertical directions. These line segments shall be denoted $b^H$ and $b^V$ respectively. By Lemma \ref{hnhd}, there is a horizontal arc $c^H$ with endpoints on $b^V$ satisfying
	\[l_q(b^V)l_q(c^H) \leq 1.\]
	Let $\gamma = c^H \cup c^V$, where ${c^V \subseteq b^V}$ is the subarc of $b^V$ connecting the endpoints of $c^H$. Note that the quantity $l_q(c^H)l_q(c^V)$ remains constant under the Teichm\"uller geodesic flow $g_t$. Therefore we may, in addition, choose $q\in\disc$ to satisfy $l_q(c^H) = l_q(c^V) \leq l_q(b^V)$. Combined with the above inequality, it follows that $l_q(c^H)$ and $l_q(c^V)$ are both at most 1 and hence
	\[l_\disc(\gamma) \leq l_q(\gamma) \leq l_q(c^H) + l_q(c^V) \leq 2.\]
	Finally, observe that $i(\beta, \gamma) \leq 1$ which implies $d_S(\beta,\gamma) \leq 2$.
	\endproof
	
	\begin{Prop}\label{nearsystole}
	 The set $V(\disc)$ is contained in the 54--neighbourhood of $\sys(\disc)$ in $\C(S)$.
	\end{Prop}
	
	\proof 
	Let $\alpha \in V(\disc)$ be a straight vertex cycle for $\disc$.
	Recall from the proof of Lemma \ref{vertexmap} that $\alpha$ must be a vertex cycle of $\td(\alpha)$. Choose some $q \in \disc$, and let $\Gamma$ be the straightening of $\tq(\alpha)$ on $q$, i.e. the embedded graph whose edges are exactly the saddle connections used by $\alpha^q$. Note that $\alpha^q$ runs over each edge of $\Gamma$ at most twice. Using an Euler characteristic argument with the Gauss-Bonnet Theorem, one can show that $\Gamma$ has at most $18(g-1)$ edges. Applying Lemma \ref{nearbicorn}, there exists some cylinder curve or $\disc$--bicorn $\beta$ satisfying $i(\alpha,\beta) \leq 36(g-1)$, and so by Lemma \ref{distbound} we have $d_S(\alpha,\beta) \leq 15$.
	By Lemma \ref{nearshort}, there is some curve $\gamma\in\C(S)$ satisfying $l_\disc(\gamma) \leq 2$ within distance 2 of $\beta$. Using Lemma \ref{shortbounded} with $L = 2$, we deduce that $\gamma$ is in the 37--neighbourhood of $\sys(\disc)$. Finally, applying the triangle inequality completes the proof.
	\endproof

	\subsection{Cylinder curves}
	
	We conclude this section by showing that the set of cylinder curves for a Teichm\"uller disc $\disc$ agrees with the systole set up to uniformly bounded Hausdorff distance in the curve graph. 
	The key ingredient for our proof is a stronger version of Bowditch's wide annulus result (Proposition \ref{wideannulus}) for flat surfaces due to Vorobets. 
				
	\begin{Prop}[\cite{Vorobets-cyl}]\label{widecylinder}
	Let $q$ be unit-area half-translation surface of genus $g \geq 2$. Then there is a flat cylinder on $q$ with width at least
	 ${\widthV = (4 \sqrt{2} (g-1) 2^{2^{32(g-1)}})^{-1}}$. \halmos
	\end{Prop}
	
	Vorobets' statement is originally for translation surfaces, but it can be generalised to half-translation surfaces by taking a branched double cover the singular points.
	
	Let $\cyl(\disc)$ and $\hcyl(\disc)$ respectively denote the set of cylinder curves and constant direction curves for any (hence all) $q\in\disc$. 
		
	\begin{Lem}\label{nearcyl}
	For any Teichm\"uller disc $\disc\subset \Teich(S_g)$, we have
	\[\hcyl(\disc) \approx_1 \cyl(\disc) \approx_{\hdSC} \sys(\disc)\]
	where $\hdSC \asymp \frac{2^{32g}}{\log g}$.
	\end{Lem}

	\proof
	We will proceed by showing the following chain of coarse inclusions:
	 \[\cyl(\disc) \subseteq \hcyl(\disc) \subseteq_1 \cyl(\disc) \subseteq_1 \sys(\disc) \subseteq_{\hdSC} \cyl(\disc)\]
	 where $\hdSC$ shall be determined later.
 
	The first inclusion holds since cylinders have constant slope. For the second inclusion, suppose $\alpha$ has constant direction (which we may assume to be vertical) on some $q\in\disc$. By shrinking the flat structure in the vertical and expanding in the horizontal direction, we can make $\alpha$ arbitrarily short and hence disjoint from a flat cylinder of width at least $\widthV$.

	If $\alpha$ is a (vertical) cylinder curve, then we can make this cylinder arbitrarily wide by expanding in the horizontal and shrinking in the vertical direction. Since systoles have length at most $\frac{2}{\sqrt{\pi}}$, it follows that $\alpha$ is disjoint from a systole on some $q \in \disc$, giving us the third inclusion.

	Finally, suppose $\alpha\in\sys(q)$ for some $q\in\disc$, and let $\beta$ be a cylinder with width at least $\widthV$ on $q$. Then $l_q(\beta) \leq \frac{1}{\widthV}$ and so by Lemma \ref{shortbounded}
	\[d_S(\alpha,\beta) \prec_{2,\diamL} \log_\Gen \left(\frac{2}{\sqrt{\pi}}\right) + \log_\Gen \left(\frac{1}{\widthV}\right) 
	\asymp \frac{2^{32g}}{\log g}\]
	which gives the fourth inclusion. One can choose $\hdSC \asymp \frac{2^{32g}}{\log g}$.
	\endproof
	
	Our bound $\hdSC$ depends on the genus of $S$, though we do not expect this to be sharp.
	This motivates the following.
	
	\begin{Que}
	What are the optimal asymptotics for the Hausdorff distance between $\cyl(\disc)$ and $\sys(\disc)$ in terms of genus? Does there exist a universal bound?
	\end{Que}
	
	Similar questions can also be posed when restricting to certain classes of half-translation surfaces, such as square-tiled surfaces, Veech surfaces, or completely periodic surfaces.

	\section{Auxiliary polygons}\label{secpolygon}

	To a curve $\alpha\in\C(S)$ and a half-translation surface $q \in \QD(S)$, we associate a convex Euclidean polygon $P_q(\alpha) \subset \R^2$ called its \emph{auxiliary polygon}. Geometric properties of $\alpha^q$ under $SL(2,\R)$--deformations of $q$ can be observed by performing the same deformations on $P_q(\alpha)$. This allows us to simplify many of our arguments in Section \ref{secbalance}. The construction works equally well if $\alpha$ is a multicurve on $S$, or a multi-arc on $S'(\disc)$, but we shall focus only on curves to simplify the exposition.
	
	Consider a geodesic representative $\alpha^q$ of $\alpha$ on $q$ which is a concatenation of saddle connections. Suppose the saddle connections $e_1, \ldots, e_k$ of $\alpha^q$ appear with multiplicities $w_1, \ldots, w_k$. Each saddle connection $e_i$ has a well-defined length and direction (up to rotation by $\pi$) on $q$. To each $e_i$ we associate a parallel vector $\uvec_i\in\R^2$ of magnitude $w_i \times \length(e_i)$. Without loss of generality, we assume the directions of the vectors $\uvec_1, \ldots, \uvec_k, -\uvec_1, \ldots, -\uvec_k$ appear in increasing anticlockwise order. (In general, this order will not agree with the order in which the saddle connections appear on $\alpha^q$.)

	\begin{Dfn}
	The \emph{auxiliary polygon} of $\alpha$ with respect to $q$ is
	\[P_q(\alpha) = \left\{\sum_i t_i \uvec_i ~|~ 0 \leq t_i \leq 1 \right\} \subset \R^2.\]
	\end{Dfn}

	Observe that $P_q(\alpha)$ is precisely the convex hull of the finite set of points $\left\{\sum_i \epsilon_i \uvec_i ~|~ \epsilon_i \in \{0, 1\} \right\}$ in $\R^2$. Thus this construction produces a convex Euclidean polygon unless all the $\uvec_i$ are parallel, in which case the polygon $P_q(\alpha)$ degenerates to a straight line segment. This situation occurs exactly when all saddle connections of $\alpha^q$ are parallel, i.e. when $\alpha^q$ has constant direction on $q$. We call $P_q(\alpha)$ \emph{degenerate} when this happens. Furthermore, if $\alpha$ is a cylinder curve then $P_q(\alpha)$ is a line segment whose direction and length is that of any geodesic representative of $\alpha$ on $q$.

	There is some ambiguity arising from the $2k$ possible cyclic permutations of the $\pm\uvec_i$, however all such choices will produce an identical polygon up to translation in $\R^2$. Indeed, $P_q(\alpha)$ is the image of the unit cube $[0,1]^k \subset \R^k$ under the linear map $f : \R^k \rightarrow \R^2$ defined by setting $f(\mathbf{e}_i) = \uvec_i$, where $\mathbf{e}_i$ for $1\leq i\leq k$ form the standard basis for $\R^k$.

	If all of the saddle connections of $\alpha^q$ are non-parallel, then as one follows $\partial P_q(\alpha)$ in an anticlockwise direction, its edges (viewed as oriented line segments) coincide exactly with $\uvec_1, \ldots, \uvec_k, -\uvec_1, \ldots, -\uvec_k$ up to cyclic permutation. In general, an edge of $\partial P_q(\alpha)$ coincides with the sum of consecutive parallel vectors. In the case where $P_q(\alpha)$ degenerates to a single straight line segment, we view $\partial P_q(\alpha)$ as a closed path which traverses the line segment once in each direction.  
	
	We may also think of $P_q(\alpha)$ as being constructed as follows: take two copies of each saddle connection of $\alpha^q$ on $q$ (counting multiplicity) and place them on the Euclidean plane (maintaining their directions) as line segments with opposing orientations. Then translate the oriented line segments in $\R^2$ so that they are positioned head to tail in order of increasing direction. This yields a closed polygonal path forming the boundary of $P_q(\alpha)$.

	For any Euclidean polygon $P \subseteq \R^2$, let $\area(P)$ denote its Euclidean area. We also define
	\[\width(P) = \sup \{|x_1 - x_2| ~|~ (x_1,y_1), (x_2,y_2) \in P\}\]
	and
	\[\height(P) = \sup \{|y_1 - y_2| ~|~ (x_1,y_1), (x_2,y_2) \in P\}.\]

	We may pull back several notions of length on $\R^2$ to $S$ via $q$: let $l_q^H(\alpha)$, $l_q^V(\alpha)$ and $l_q(\alpha)$ respectively denote the length of $\alpha^q$ with respect to $|dx|$, $|dy|$ and the Euclidean metric. The following is immediate.

	\begin{Lem}\label{polydims}
	Let $q$, $\alpha$ and $P_q(\alpha)$ be as above. Then
		\begin{itemize}
			\item	$\width(P_q(\alpha)) = l_q^H(\alpha)$, 
			\item $\height(P_q(\alpha)) = l_q^V(\alpha)$, and 
						\item$l(\partial P_q(\alpha)) = 2\,l_q(\alpha)$. \halmos
		\end{itemize}
	\end{Lem}

	\begin{Lem}
	For any $A \in SL(2,\R)$, the polygons $A \cdot P_q(\alpha)$ and $P_{A \cdot q}(\alpha)$ agree up to translation in $\R^2$. In particular, we have
	\[\area (P_{A \cdot q}(\alpha)) = \area (P_q(\alpha)).\]
	\end{Lem}
	
	\proof
	Observe that elements of $SL(2,\R)$ preserve the cyclic order on the directions of saddle connections $e_i$, and hence the cyclic order of the $\pm\uvec_i$. Then, up to translation, we have 
	\[A \cdot P_q(\alpha) = A \cdot \left\{\sum_i t_i \uvec_i ~|~ 0 \leq t_i \leq 1 \right\} = \left\{\sum_i t_i A \uvec_i ~|~ 0 \leq t_i \leq 1 \right\} = P_{A \cdot q}(\alpha). \]
	The second claim follows from the fact that Euclidean area is preserved under $SL(2,\R)$--deformations.
	\endproof

	For any curve $\alpha\in\C(S)$ and Teichm\"uller disc $\disc$, define its \emph{polygonal area} with respect to $\disc$ to be
	\[\area_\disc(\alpha) = \area (P_q(\alpha))\]
	 where $q$ is any unit area half-translation surface in $\disc$. We also define the \emph{infimal length}
	\[l_\disc(\alpha) = \inf_{q\in\disc}\, l_q(\alpha)\]
	of $\alpha$ with respect to $\disc$.
	
	\begin{Prop}\label{arealength}
	With $\alpha$ and $\disc$ as above, we have
	\[\pi\,\area_\disc(\alpha) \leq (l_\disc(\alpha))^2 \leq 8\, \area_\disc(\alpha).\]
	In particular, $l_\disc(\alpha) = 0$ if and only if $\alpha \in \hcyl(\disc)$. Furthermore, the infimal length of $\alpha$ is realised if and only if it is positive.
	\end{Prop}
	
	Thus, one can estimate the infimal length of a curve $\alpha$ over $\disc$ by computing the area of $P_q(\alpha)$ for any $q \in \disc$.

	To prove this result, we require the following \emph{Round Polygon Lemma}.
	Let $B^1(r)$ and $B^\infty(r)$ respectively denote the balls of radius $r$ about the origin in $\R^2$ with respect to the $l^1$ and $l^\infty$ norms.
		
	\begin{Lem}\label{round}
	Suppose $P$ is a convex (non-degenerate) Euclidean polygon in $\R^2$ with $\pi$--rotational symmetry about the origin. Then there exists $A \in SL(2,\R)$ such that
	\[B^1(r) \subseteq A\cdot P \subseteq B^\infty(r) \]
	for some $r>0$.
	\end{Lem}
	
	\begin{figure}[h]
	    \centering{
	    \resizebox{75mm}{!}{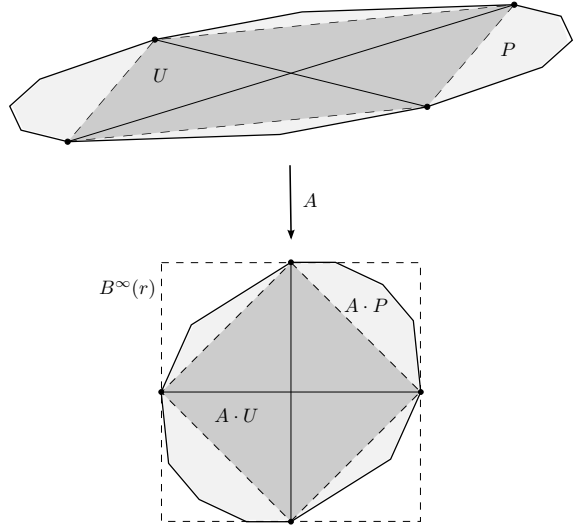}
	    \caption{Applying $A \in \SL(2,\R)$ to the polygon $P$ to make it ``round''.}
	    \label{fig:roundpolygon}
	    }
	\end{figure}
	
	\proof
	Consider the set of parallelograms spanned by a pair of diagonals of $P$ which pass through the origin, and choose one such parallelogram $U \subseteq P$ with largest area (see Fig. \ref{fig:roundpolygon}). Deform $P$ using some $A \in SL(2,\R)$ so that $A \cdot U = B^1(r)$ for some $r>0$. We claim that $A\cdot P \subseteq B^\infty(r)$. Supposing otherwise, $A\cdot P$ must have a corner $(x,y)\in\R^2$ lying outside $B^\infty(r)$. Without loss of generality, we may assume $|x| > r$. Since $A \cdot P$ has $\pi$--rotational symmetry about its centre, the point $(-x,-y)$ is also a corner of $A \cdot P$. Let $U'\subseteq A\cdot P$ be the parallelogram spanned by the points $\pm (x,y)$ and $\pm (0,r)$, the latter pair being corners of $B^1(r)$ and hence $A \cdot P$. Now
	\[\area(A^{-1} \cdot U') = \area(U') = 2r|x| > 2r^2 = \area(A\cdot U) = \area(U), \]
	contradicting the maximality assumption on $U$.
	\endproof

	\proofof{Proposition \ref{arealength}}
	To prove the lower bound on $l_\disc(\alpha)$, we apply the well-known isoperimetric inequality for the plane to $P_q(\alpha)$ for every $q \in \disc$: A planar region $U \subset \R^2$ enclosed by an embedded curve of length $L = 2l_q(\alpha)$ must satisfy $4\pi\area(U) \leq L^2$.
	
	For the other direction, we may apply Lemma \ref{round} to find some $q\in\disc$ so that $B^1(r) \subseteq P_q(\alpha) \subseteq B^\infty(r)$ for some $r>0$. We divide the boundary of $P = P_q(\alpha)$ into four subpaths connecting adjacent corners of $B^1(r)$. The subpath $\eta$ in the first quadrant is a concatenation of straight line segments connecting consecutive points with co-ordinates 
	$(r,0) = (x_0, y_0), (x_1, y_1), \ldots, (x_m, y_m) = (0, r)$.
	Since $P$ is convex, it follows that $x_0 \geq x_1 \geq \ldots \geq x_m $ and $y_0 \leq y_1 \leq \ldots \leq y_m$, and hence
	\[l(\eta) \leq l^1(\eta) = \sum_{i=1}^m |x_i - x_{i-1}| + |y_i - y_{i-1}| = 2r.\]
	Applying this to the other subpaths, we deduce
	$2l_\disc(\alpha) = l(\partial P) \leq 8r$. We finally combine this with the inequality $\area(P) \geq \area(B^1(r)) = 2r^2$ to obtain the desired upper bound.
	
	Observe that $l_\disc(\alpha) = 0$ if and only if $\area_\disc(\alpha) = 0$, which occurs precisely when $P_q(\alpha)$ is degenerate for any (hence all) $q \in \disc$.
	\halmos
	
	\section{Balance points on Teichm\"uller discs}\label{secbalance}
	
	In this section, we generalise Masur and Minsky's notion of balance time on Teichm\"uller geodesics to Teichm\"uller discs. In particular, we prove the existence of a \emph{balance point} which coarsely determines the balance time for \emph{every} Teichm\"uller geodesic contained inside a common Teichm\"uller disc. The auxiliary polygon from Section \ref{secpolygon} plays a key role in our proofs.
	
	Let $\disc$ be a Teichm\"uller disc and $\alpha\in\C(S)$ be a curve. 
	If $\alpha\in \hcyl(\disc)$, we define the \emph{balance point} of $\alpha$ with respect to $\disc$ to be the projectivised measured foliation $\F \in \partial \disc \subset \PMF(S)$ in the direction of $\alpha$.
	If $\alpha\not\in\hcyl(\disc)$, we say $q\in\disc$ is a balance point of $\alpha$ with respect to $\disc$ if the auxiliary polygon $P_q(\alpha)$ is round in the sense of Lemma \ref{round}:
	there is some $r>0$ such that $B^1(r) \subseteq P_q(\alpha) \subseteq B^\infty(r)$.
	Balance points always exist by Lemma \ref{round}. Write $\G_\alpha$ for the balance time of $\alpha$ along a Teichm\"uller geodesic $\G$. If $\alpha$ is completely horizontal or completely vertical with respect to $\G$, then we set $\G_\alpha$ to be the endpoint of $\G$ corresponding to the direction of $\alpha$.

	\begin{Prop}\label{balance}
	 Let $X\in\disc\cup \partial\disc$ be a balance point for a curve $\alpha\in\C(S)$ on a Teichm\"uller disc $\disc$. 
		 	 	 	 For any Teichm\"uller geodesic $\G \subset \disc$, let $p\in \G$ be the nearest point projection of $X$ to $\G$ in $\disc$. Then
	  $d_\disc(p,\G_\alpha) \leq \log 2$.
	 In particular, if $\G$ passes through $X$ then $d_\disc(X,\G_\alpha) \leq \log 2$.
		\end{Prop}
	
	It follows that balance points are coarsely unique: If $X$ and $X'$ are balance points for $\alpha$ on $\disc$, then $d_\disc(X,X') \leq 2 \log 2$. If a balance point is in $\partial \disc$, then it is unique by definition. Combining the above with Lemma \ref{lengthconvex}, we deduce the following.
	
	\begin{Cor}
	 Given a Teichm\"uller disc $\disc$ and a curve $\alpha\in\C(S)$, suppose $l_\disc(\alpha)$ is realised at ${m \in \disc \cup \partial \disc}$. Then 
	  $d_\disc(X,m) \leq \cosh^{-1}2 + \log 2$, where $X$ is the balance point of $\alpha$ on $\disc$. \halmos
	\end{Cor}

	\begin{Thm}\label{balancesystole}
	 Let $X$ be a balance point for a curve $\alpha\in\C(S)$ on a Teichm\"uller disc $\disc$. Then any systole on $X$ is universally close to any nearest point projection of $\alpha$ to $\sys(\disc)$ in $\C(S)$.
	\end{Thm}

	\subsection{Proof of Proposition \ref{balance}}

	For this section, fix a Teichm\"uller disc $\disc$ and a curve $\alpha\in\C(S)$. We will deal with two separate cases, corresponding to whether the geodesic representative $\alpha^q$ on any (hence all) $q\in\disc$ has constant direction. We shall be keeping track of directions on $q$.
	
		Suppose $\alpha \in \hcyl(\disc)$. The balance point of $\alpha$ is the foliation $\F\in\partial\disc$ which is in the same direction as $\alpha$ on every $q\in\disc$. Any Teichm\"uller geodesic $\G\subset\disc$ with an endpoint at $\F\in\partial\disc$ has $\F$ as its vertical foliation (once $\G$ is oriented so that the forwards direction points towards $\F$). Since $\alpha$ is completely vertical along $\G$ then, by definition, $\F$ is also the balance point of $\alpha$ on $\G$.	
	Now let $\G\in\disc$ be any Teichm\"uller geodesic with neither endpoint at $\F$. The nearest point projection $p\in\G$ of $\F$ to $\G$ is the endpoint of the unique ray from $\F$ which meets $\G$ at a right angle on $\disc$. Therefore, $\F$ makes an angle of $\frac{\pi}{4}$ with both the horizontal and vertical directions on $p$. Since $\alpha$ is in the direction of $\F$, it follows that $l_p^H(\alpha) = \frac{1}{\sqrt{2}}l_p(\alpha) = l_p^V(\alpha)$ and thus $\G_\alpha = p$.

		Now assume $\alpha$ has non-constant direction on $\disc$ and let $q\in\disc$ be a balance point. 
	Let $\G \subset \disc$ be a Teichm\"uller geodesic. Apply a rotation $\rho_\theta \in SO(2,\R)$ to $q$ so that the foliations corresponding to $\G$ make angles of $\pm \phi$ with the vertical, for some $0 < \phi \leq \frac{\pi}{4}$. 
	Using the Teichm\"uller geodesic flow to respectively stretch and shrink the horizontal and vertical directions on $\rho_\theta \cdot q$ by a factor of $\sqrt{\cot\phi}$, we obtain a new half-translation surface $p\in \G$ on which the foliations corresponding to $\G$ are perpendicular. Moreover, these foliations make an angle of $\frac{\pi}{4}$ with the horizontal and vertical. Therefore the Teichm\"uller geodesic in $\disc$ connecting $q$ to $p$ is perpendicular to $\G$, and so $p$ is the nearest point projection of $q$ to $\G$ in $\disc$.
	Let $l^H(\alpha)$ and $l^V(\alpha)$ be the horizontal and vertical lengths of $\alpha$ on $\rho_{\frac{\pi}{4}}\cdot p$. These are exactly the lengths of $\alpha$ on $p$ with respect to the foliations corresponding to $\G$. Thus, to prove $p$ is near the balance point of $\alpha$ to $\G$ we must show that these lengths are almost equal.

			\begin{center}
	\begin{figure}[h]
	    \centering{
	    \resizebox{140mm}{!}
	     {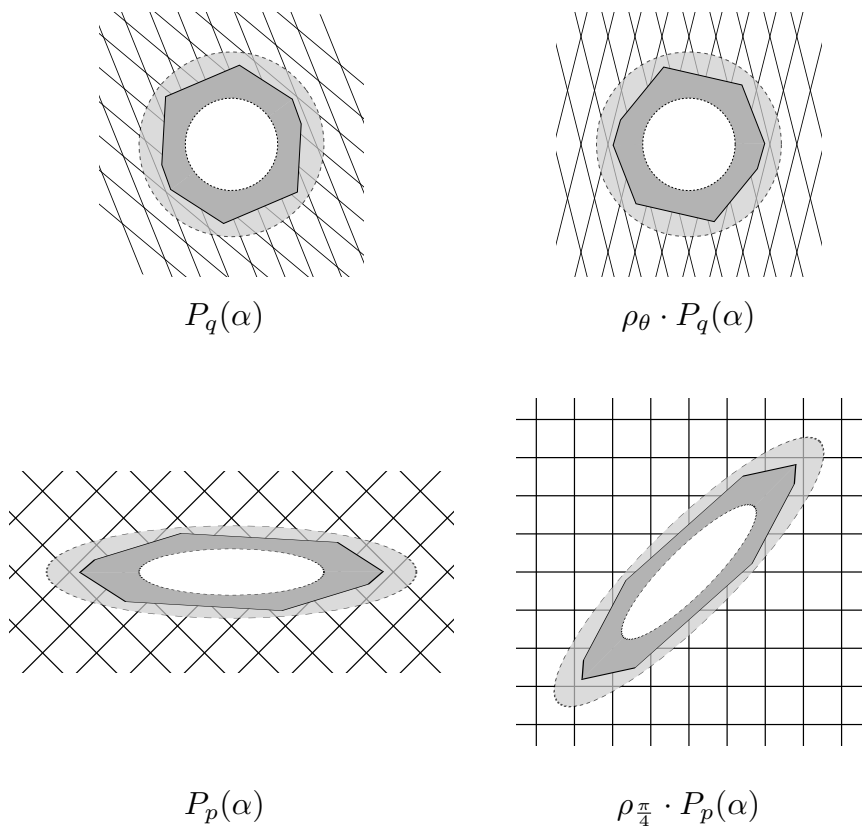}
	    \caption{The auxiliary polygons for $\alpha$, corresponding to the half-translation structures $q$, $\rho_\theta \cdot q$, $p$, and $\rho_{\frac{\pi}{4}}\cdot p$ respectively, are nested between pairs of ellipses. The grid lines indicate the pair of transverse directions corresponding to $\G$.}
	    \label{fig:ellipse}
	    }
	\end{figure}
	\end{center}

	Our strategy is to apply appropriate $\SL(2,\R)$--deformations to $P = P_q(\alpha)$, as illustrated in Figure \ref{fig:ellipse}, in order to obtain estimates for $l^H(\alpha)$ and $l^V(\alpha)$ via Lemma \ref{polydims}. By Proposition \ref{arealength}, the auxiliary polygon $P = P_q(\alpha)$ is non-degenerate and satisfies $B^1(r) \subseteq P \subseteq B^\infty(r)$ for some $r>0$. Observe that $\rho_\theta\cdot P$ can be nested between a pair of concentric Euclidean circles
	\[B^2\left(R\right) \subseteq \rho_\theta \cdot B^1(r) \subseteq \rho_\theta\cdot P \subseteq \rho_\theta \cdot B^\infty(r) \subseteq  B^2(2R),\]
	where $R = \frac{r}{\sqrt{2}}.$ The polygon $P_p(\alpha)$ can be obtained by respectively stretching and shrinking $\rho_\theta\cdot P$ by a factor of $\sqrt{\cot\phi}$ in the horizontal and vertical directions. By applying the same deformations to the circles, $P_p(\alpha)$ can be nested between a pair of ellipses
	\[ E_{inner} \subseteq P_p(\alpha) \subseteq E_{outer}\]
	defined by the equations
	\[\frac{x^2}{\cot\phi} + (\cot\phi) y^2 \leq R^2 \qquad\textrm{ and }\qquad \frac{x^2}{\cot\phi} + (\cot\phi) y^2 \leq 4R^2\]
	respectively.

	We then rotate $P_p(\alpha)$ through an angle of $\frac{\pi}{4}$ to obtain $\rho_{\frac{\pi}{4}}\cdot P_p(\alpha)$. By Lemma \ref{polydims}, we know that
	\[l^H(\alpha) = \width(\rho_{\frac{\pi}{4}}\cdot P_p(\alpha)) \qquad\textrm{ and }\qquad l^V(\alpha) = \height(\rho_{\frac{\pi}{4}}\cdot P_p(\alpha)).\]
	Rotating the ellipses $E_{inner}$ and $E_{outer}$ through an angle of  $\frac{\pi}{4}$ yields the following estimates:
	\[\width(\rho_{\frac{\pi}{4}}\cdot E_{inner}) \leq l^H(\alpha) \leq \width(\rho_{\frac{\pi}{4}}\cdot E_{outer}) = 2\,\width(\rho_{\frac{\pi}{4}}\cdot E_{inner})\]
	and
	\[\height(\rho_{\frac{\pi}{4}}\cdot E_{inner}) \leq l^V(\alpha) \leq \height(\rho_{\frac{\pi}{4}}\cdot E_{outer}) = 2\,\height(\rho_{\frac{\pi}{4}}\cdot E_{inner}).\]	
	Note that the ellipses $\rho_{\frac{\pi}{4}}\cdot E_{inner}$ and $\rho_{\frac{\pi}{4}}\cdot E_{outer}$ have reflective symmetry about the line $y=x$. This implies
	\[\width(\rho_{\frac{\pi}{4}}\cdot E_{inner}) = \height(\rho_{\frac{\pi}{4}}\cdot E_{inner}).\]
	Combining this with the preceding inequalities yields 
	\[\frac{1}{2}l^H(\alpha) \leq l^V(\alpha) \leq 2 l^H(\alpha),\]
	and therefore $p$ is within a distance of $\log 2$ of the balance time of $\alpha$ to $\G$.
	\halmos
	
	\subsection{Projecting to the systole set}
	
	In this section, we prove Theorem \ref{balancesystole} which generalises the following theorem for Teichm\"uller geodesics to Teichm\"uller discs. Let $\G_\alpha$ denote the balance time of a curve $\alpha$ to a Teichm\"uller geodesic $\G$.
	
	\begin{Thm}[\cite{MM1}]\label{balproj}
	For any Teichm\"uller geodesic $\G\subset\Teich(S)$, the systole set $\sys(\G)$ is a uniform reparameterised {quasigeodesic} in $\C(S)$.
	Furthermore, the relation $\alpha \mapsto \sys(\G_\alpha)$ is a uniformly coarse Lipschitz retract from $\C(S)$ to $\sys(\G)$. \halmos
	\end{Thm}

	In \cite{bhb-unif}, Bowditch gives universal bounds on the quasiconvexity and Lipschitz constants in the case of Teichm\"uller geodesics $\G(\beta,\beta')$ arising from a half-translation structure dual to a filling pair of weighted multicurves $\beta$ and $\beta'$. He defines a ``line'' between pairs of weighted multicurves $\beta$ and $\beta'$ in $\C(S)$ in terms of intersection numbers as follows \cite{bhb-int}.
	Set 
	\[\L_t(\beta,\beta') = \{\gamma \in \C(S) ~|~  e^ti(\beta,\gamma) + e^{-t}i(\beta',\gamma) \leq \sL_0 \sqrt{i(\beta,\beta')} \} \]
	and $\L(\beta',\beta) = \bigcup_{t\in \R} \L_t(\beta,\beta')$,
	where $\sL_0$ is a suitable universal constant. Let $q$ be the half-translation structure dual to $\beta$ and $\beta'$, where its horizontal and vertical directions respectively agree with $\beta$ and $\beta'$. Then $\L_t(\beta, \beta')$ is precisely the set of bounded length curves on $g_t \cdot q$ under the $L_1$--metric. Applying Lemma \ref{shortbounded}, we see that $\L_t(\beta,\beta') \approx \sys(g_t\cdot q)$ and hence $\L(\beta,\beta') \approx \sys(\G(\beta,\beta'))$.
	
	Bowditch proves that this system of lines satisfies a ``slim triangles'' condition with universal constants which, by a criterion of Masur and Schleimer \cite{MS-disc}, implies the curve graph is universally hyperbolic and that the line $\L(\beta,\beta')$ is a universal reparametrised  quasigeodesic between $\beta$ and $\beta'$. 
	Furthermore, he approximates ``coarse centres'' of geodesic triangles using balance times. We can reinterpret this in terms of nearest point projections: If $\G_\alpha$ is the balance time of $\alpha$ along $\G = \G(\beta,\beta')$, then $\sys(\G_\alpha)$ is universally close to any nearest point projection of $\alpha$ to $\sys(\G)$ in $\C(S)$.
		
	We now outline a proof extending Bowditch's results to arbitrary Teichm\"uller geodesics via a limiting argument.
	
	\begin{Prop}\label{sysuniv}
	 For any Teichm\"uller geodesic $\G \subset \Teich(S)$, the systole set $\sys(\G)$ is $\qcS$--quasiconvex for some universal constant $\qcS$. Furthermore, the operation $\alpha \mapsto \sys(\G_\alpha)$ agrees with the nearest point projection from $\C(S)$ to $\sys(\G)$ up to a universal error $\projS$.
	\end{Prop}

	\proof
	Suppose $\G$ is an arbitrary Teichm\"uller geodesic whose endpoints correspond to a pair of transverse measured foliations $\lambda$ and $\lambda'$. Define $\L_t(\lambda, \lambda')$ and $\L(\lambda, \lambda')$ using the same intersection number conditions as above. Using a suitable parameterisation, these sets respectively agree with $\sys(\G_t)$ and $\sys(\G)$ up to universal Hausdorff distance. Let $\beta_n$ and $\beta'_n$ be sequences of weighted multicurves converging in the space of measured foliations $\MF(S)$ to $\lambda$ and $\lambda'$ respectively. Appealing to continuity of intersection number on $\MF(S) \times \MF(S)$ \cite{Bonahon-currents} and Lemma \ref{shortbounded}, for each $t \in \R$ we have $\L_t(\beta_n, \beta'_n) \approx \L_t(\lambda, \lambda')$ for all $n$ sufficiently large. 
	 Since each $\L(\beta_n, \beta'_n)$ is universally quasiconvex \cite{bhb-unif}, it follows that $\L(\lambda, \lambda')$ and hence $\sys(\G)$ are also universally quasiconvex.
	
	To prove the second claim, one can again appeal to continuity of intersection number to show that $\sys(\G(\beta_n,\beta'_n)_\alpha)$ will eventually agree with $\sys(\G_\alpha)$ up to universal Hausdorff distance for all sufficiently large $n$. Applying Bowditch's nearest point projection result completes the proof.
	\endproof
	
	Using the above results, we can give an alternative proof that systole sets are universally quasiconvex, though with weaker effective control over the constants.
	
	\begin{Cor}
	 For any Teichm\"uller disc $\disc \subset \Teich(S)$, the systole set $\sys(\disc)$ is $\qcS$--quasiconvex.
	\end{Cor}
	
	\proof
	Given any pair of curves $\alpha, \beta \in \sys(\disc)$, let $X,Y \in \disc$ be points such that $\alpha \in \sys(X)$ and $\beta \in \sys(Y)$. Let $\G\subset\disc$ be a Teichm\"uller geodesic connecting $X$ and $Y$. Then any geodesic in $\C(S)$ connecting $\alpha$ to $\beta$ must lie within a distance $\qcS$ of $\sys(\G) \subseteq \sys(\disc)$.
	\endproof

	\proofof{Theorem \ref{balancesystole}}
	Let $\alpha\in\C(S)$ be a curve and $X$ a balance point with respect to a Teichm\"uller disc $\disc$.
	Let $\gamma$ be a nearest point projection of $\alpha$ to $\sys(\disc)$ in $\C(S)$. Then $\gamma$ is a systole for some $Y\in\disc$. If $Y$ coincides with $X$ then we are done, so suppose otherwise. Let $\G$ be an infinite Teichm\"uller geodesic passing through $X$ and $Y$. By Theorem \ref{balproj} and Proposition \ref{sysuniv}, $\sys(\G_\alpha)$ is universally close to $\gamma$. By Theorem \ref{balance}, we have $d_\disc(X,\G_\alpha) \leq  \log 2$ and so by Lemma \ref{syslip}, it follows that $\sys(X)$ is uniformly close to $\sys(\G_\alpha)$.
	\halmos

	The above proof also works if we let $X$ be the point which minimises the length of $\alpha$ on $\disc$.

	\section{Balance points and curve decompositions}\label{secvertexbalance}

	Fix a Teichm\"uller disc $\disc \subset \Teich(S)$. Let $\tau\in\TT(\disc)$ be a straight train track with respect to $\disc$, and suppose $\alpha$ is a multicurve carried by $\tau$. By Theorem \ref{vertexdecomp}, we can write $w_\alpha = \frac{1}{2}\sum_{v \in V(\tau)} m_v w_v$ where the $m_v$ are non-negative integers.
	Let $X$ and $X_v$ respectively be balance points for $\alpha$ and $v$ on $\disc$. Let $H = H(\tau)$ be the convex hull of the $X_v$'s in $\disc$. 
	
	The goal of this section is to show the following connection between balance points and straight vertex cycle decompositions.
	
	\begin{Prop}\label{hullvertex}
	The sets $\sys(X)$, $\sys(H)$, and $V(\tau)$ agree up to universally bounded Hausdorff distance in $\C(S)$. In particular, $\sys(H)$ has universally bounded diameter.
	\end{Prop}
	
	Combining this result with Corollary \ref{sysnearvertex}, Proposition \ref{nearsystole}, and Theorem \ref{balancesystole}, we can give an alternative proof to Proposition \ref{vertexproj}, albeit with weaker control over the constants.
	
	\begin{Lem}
	Let $\G$ be a Teichm\"uller geodesic in $\disc$. Let $q$ and $q_v$ be the respective balance times of $\alpha$ and $v$ on $\G$, for all $v \in V(\tau)$. If $I \subseteq \G$ is the minimal subinterval which contains all the $q_v$'s, then $q$ is also contained in $I$.
	\end{Lem}

	\proof
	We may assume the horizontal directions on $q$ and all the $q_v$'s are the same, i.e. they correspond to the same endpoint of $\G$.
	Suppose the claim is false. Then the $q_v$'s must all lie on the same component of $\G - \{q\}$. This means that either $l_q^H(v) > l_q^V(v)$ for all $v$; or $l_q^H(v) < l_q^V(v)$ for all $v$. It follows that either \[l_q^H(\alpha) = \frac{1}{2}\sum_v m_v l_q^H(v) > \frac{1}{2}\sum_v m_v l_q^V(v) = l_q^V(\alpha);\]
	or $l_q^H(\alpha) < l_q^V(\alpha)$. Then $\alpha$ is not balanced at $q$, a contradiction.
	\endproof

	\begin{Lem}\label{nearhull}
	Given $\tau \in \TT(\disc)$, let $H = H(\tau)$ be as above. Then for any multicurve $\alpha$ carried by $\tau$, any balance point $X$ of $\alpha$ on $\disc$ is in the $2\log 2$--neighbourhood of $H$. 	\end{Lem}

	\proof
	If $X$ is in $H$ then we are done, so suppose otherwise. Let $Y$ be the unique closest point projection of $X$ to $H$, and $\G$ be the infinite Teichm\"uller passing through $X$ and $Y$. Let $q$, $q_v$ and $I$ be as in the previous lemma for $\G$. By Proposition \ref{balance}(1), 
	\[d_\disc(X,H) = d_\disc(X,Y) \leq d_\disc(X,q) + d_\disc(q,Y) \leq \log 2 + d_\disc(q,Y)\] 	and so it suffices to prove $d_\disc(q,Y) \leq \log 2$.
	
	Let $J$ and $J'$ be the two closed rays contained in $\G$ which have $Y$ as their endpoint. Assume $X\in J'$. Let $\pi_\G : \Delta \rightarrow \G$ be the closest point projection map, and let $p_v = \pi_\G(X_v)$. Using elementary hyperbolic geometry, one can show that $\pi_\G(H) \subseteq J$. By Proposition \ref{balance}(2), we have $d_\disc(q_v, p_v) \leq \log 2$ which implies $I \subseteq_{\log 2} J$. Applying the previous lemma, we have $q \in I$ and so
	\[d_\disc(q, J) \leq \log 2.\]	
	On the other hand,
	\[d_\disc(q, J') \leq d_\disc(q,X) \leq \log 2\] 	by Proposition \ref{balance}(1). Combining these two bounds, we deduce that $q$ is within distance $\log 2$ of the common endpoint of $J$ and $J'$, namely $Y$.	
	\endproof
	
	\begin{Lem}\label{vertexsystole}
	Let $v \in V(\disc)$ be a straight vertex cycle for $\disc$ and let $X_v \in \Delta$ be a balance point of $v$. Then $v$ is universally close to $\sys(X_v)$ in $\C(S)$.
	\end{Lem}
	
	\proof
		Let $\gamma$ be any nearest point projection of $v$ to $\sys(\disc)$. By Proposition \ref{nearsystole} and Theorem \ref{balancesystole} respectively, we deduce
	\[d_S(v, \gamma) = d_S(v, \sys(\Delta)) \leq 54\]
	and
	\[d_S(\gamma, \sys(X_v)) \leq \projS.\]
	The result follows by applying the triangle inequality.
	\endproof

	\proofof{Proposition \ref{hullvertex}}
	By Lemmas \ref{syslip} and \ref{nearhull}, $\sys(X)$ lies in a universally bounded neighbourhood of $\sys(H)$.
	Recall from Theorem \ref{vertexdiam} that $V(\tau)$ has diameter at most $\lipV \leq 14$ in $\C(S)$. Therefore, to universally bound the pairwise Hausdorff distances between $\sys(X)$, $\sys(H)$ and $V(\tau)$ in $\C(S)$, it suffices to show that $\sys(H)$ is contained in a universally bounded neighbourhood of $V(\tau)$.

	Recall that $\disc$ is isometric to the hyperbolic plane with curvature $-4$. Thus $H$ is contained in a universally bounded neighbourhood of the union of all Teichm\"uller geodesics connecting $X_v$ and $X_{v'}$, for all pairs $v,v'\in V(\tau)$. By Lemma \ref{syslip}, Theorem \ref{balproj}, and Proposition \ref{sysuniv}, the relation $\sys : \disc \rightarrow \C(S)$ is coarsely Lipschitz and sends Teichm\"uller geodesics to reparametrised quasigeodesics with universal quasiconvexity constants. Therefore $\sys(H)$ is contained in a universally bounded neighbourhood of the set $\{\sys(X_v) ~|~ v\in V(\tau)\}$. By the previous lemma, $\sys(X_v)$ is universally close to $v$ and we are done.
				\halmos

	\bibliography{mybib}		

\providecommand{\bysame}{\leavevmode\hbox to3em{\hrulefill}\thinspace}
\providecommand{\MR}{\relax\ifhmode\unskip\space\fi MR }
\providecommand{\MRhref}[2]{%
  \href{http://www.ams.org/mathscinet-getitem?mr=#1}{#2}
}
\providecommand{\href}[2]{#2}
\begin{thebibliography}{HPW15}

\bibitem[Aou13]{Aougab-unif}
Tarik Aougab, \emph{Uniform hyperbolicity of the graphs of curves}, Geom.
  Topol. \textbf{17} (2013), no.~5, 2855--2875. \MR{3190300}

\bibitem[Bon88]{Bonahon-currents}
Francis Bonahon, \emph{The geometry of {T}eichm\"uller space via geodesic
  currents}, Invent. Math. \textbf{92} (1988), no.~1, 139--162. \MR{931208
  (90a:32025)}

\bibitem[Bow06]{bhb-int}
Brian~H. Bowditch, \emph{Intersection numbers and the hyperbolicity of the
  curve complex}, J. Reine Angew. Math. \textbf{598} (2006), 105--129.
  \MR{2270568 (2009b:57034)}

\bibitem[Bow14]{bhb-unif}
\bysame, \emph{Uniform hyperbolicity of the curve graphs}, Pacific J. Math.
  \textbf{269} (2014), no.~2, 269--280. \MR{3238474}

\bibitem[CRS14]{CRS-unif}
Matt Clay, Kasra Rafi, and Saul Schleimer, \emph{Uniform hyperbolicity of the
  curve graph via surgery sequences}, Algebr. Geom. Topol. \textbf{14} (2014),
  no.~6, 3325--3344. \MR{3302964}

\bibitem[DLR10]{DLR-length}
Moon Duchin, Christopher~J. Leininger, and Kasra Rafi, \emph{Length spectra and
  degeneration of flat metrics}, Invent. Math. \textbf{182} (2010), no.~2,
  231--277. \MR{2729268 (2011m:57022)}

\bibitem[Ham]{ham-disc}
Ursula Hamenstaedt, \emph{Geometry of graphs of discs in a handlebody: Surgery
  and intersection}, preprint, (2014), arXiv:1101.1843v4.

\bibitem[Hem01]{Hempel-cc}
John Hempel, \emph{3-manifolds as viewed from the curve complex}, Topology
  \textbf{40} (2001), no.~3, 631--657. \MR{1838999 (2002f:57044)}

\bibitem[HPW15]{HPW-unicorn}
Sebastian Hensel, Piotr Przytycki, and Richard C.~H. Webb, \emph{1-slim
  triangles and uniform hyperbolicity for arc graphs and curve graphs}, J. Eur.
  Math. Soc. (JEMS) \textbf{17} (2015), no.~4, 755--762. \MR{3336835}

\bibitem[Mas85]{Maskit-comparison}
Bernard Maskit, \emph{Comparison of hyperbolic and extremal lengths}, Ann.
  Acad. Sci. Fenn. Ser. A I Math. \textbf{10} (1985), 381--386. \MR{802500
  (87c:30062)}

\bibitem[Min01]{Minsky-qc}
Yair~N. Minsky, \emph{Bounded geometry for {K}leinian groups}, Invent. Math.
  \textbf{146} (2001), no.~1, 143--192. \MR{1859020 (2004f:30032)}

\bibitem[MM99]{MM1}
Howard~A. Masur and Yair~N. Minsky, \emph{Geometry of the complex of curves.
  {I}. {H}yperbolicity}, Invent. Math. \textbf{138} (1999), no.~1, 103--149.
  \MR{1714338 (2000i:57027)}

\bibitem[MM00]{MM2}
\bysame, \emph{Geometry of the complex of curves. {II}. {H}ierarchical
  structure}, Geom. Funct. Anal. \textbf{10} (2000), no.~4, 902--974.
  \MR{1791145 (2001k:57020)}

\bibitem[MM04]{MM3}
\bysame, \emph{Quasiconvexity in the curve complex}, In the tradition of
  {A}hlfors and {B}ers, {III}, Contemp. Math., vol. 355, Amer. Math. Soc.,
  Providence, RI, 2004, pp.~309--320. \MR{2145071 (2006a:57022)}

\bibitem[MMS12]{MMS-train}
Howard {Masur}, Lee {Mosher}, and Saul {Schleimer}, \emph{{On train-track
  splitting sequences.}}, {Duke Math. J.} \textbf{161} (2012), no.~9,
  1613--1656.

\bibitem[Mos03]{Mosher-train}
Lee Mosher, \emph{Train track expansions of measured foliations},
  \texttt{http://andromeda.rutgers.edu/\textasciitilde mosher/}, 2003.

\bibitem[MS13]{MS-disc}
Howard Masur and Saul Schleimer, \emph{The geometry of the disk complex}, J.
  Amer. Math. Soc. \textbf{26} (2013), no.~1, 1--62. \MR{2983005}

\bibitem[PH92]{pennerharer}
R.~C. Penner and J.~L. Harer, \emph{Combinatorics of train tracks}, Annals of
  Mathematics Studies, vol. 125, Princeton University Press, Princeton, NJ,
  1992. \MR{1144770 (94b:57018)}

\bibitem[Vor]{Vorobets-cyl}
Yaroslav Vorobets, \emph{Periodic geodesics on translation surfaces}, preprint,
  (2003), arXiv:math/0307249.

\end{thebibliography}
	\bibliographystyle{amsalpha}

	\end{document}